\documentclass[10pt]{amsart}
            
\usepackage{amsmath}               
\usepackage{amsfonts}              
\usepackage{amsthm}                
\usepackage[utf8]{inputenc}
\usepackage{leftidx}
\usepackage{tikz-cd}
\usepackage{mathtools}

\newtheorem{tw}{Theorem}[section]
\newtheorem{lem}[tw]{Lemma}
\newtheorem{defi}[tw]{Definition}
\newtheorem{prop}[tw]{Proposition}
\newtheorem{rem}[tw]{Remark}
\newtheorem{cor}[tw]{Corollary}

\newtheorem{obs}[tw]{Observation}

\newcommand{\Exterior}{\mathchoice{{\textstyle\bigwedge}}%
    {{\bigwedge}}%
    {{\textstyle\wedge}}%
    {{\scriptstyle\wedge}}}

\title[A Spectral Sequence for Actions of Compact Lie Groups]{A Spectral Sequence for Equidimensional Actions of Compact Lie Groups}
\author{Paweł Raźny}
\address{Institute of Mathematics \\
	Faculty of Mathematics and Computer Science \\
	Jagiellonian University in Cracow
	}
\email{pawel.razny@uj.edu.pl}


\begin{document}
\begin{abstract}
In this article we provide a version of the Leray-Serre spectral sequence for equidimensional (i.e. smooth with all orbits of the same dimension) actions of compact connected Lie groups on compact manifolds. The main part of this article consists of the proof of the description of the second page of said spectral sequence. This description provides a link between the cohomology of the orbit space (basic cohomology of the foliation by orbits) the Lie algebra cohomology of the appropriate pair $(\mathfrak{g},\mathfrak{h})$ representing the cohomology of a generic orbit and the de Rham cohomology of the manifold. Due to the somewhat technical nature of the general description we have provided in the penultimate section a thorough study of special cases in which the sequence can be greatly simplified. In particular, vast simplifications can be obtained if the manifold $M$ on which the group acts is assumed to be simply connected or if the acting Lie group has some nice properties. In the final section, we show how to use a blow up process to use our sequence when the action is not equidimensional. We apply this method to give a topological obstruction to the existence Lie group actions on certain manifolds.
\end{abstract}
\maketitle
\section{Introduction}
The purpose of this article is to provide an anologue of the Leray-Serre spectral sequence for a certain class of Lie group actions (extending the results form \cite{My4}) connecting in a meaningful way the cohomology of the orbit space (basic cohomology of the foliation by orbits) the Lie algebra cohomology of the appropriate pair $(\mathfrak{g},\mathfrak{h})$ representing the cohomology of a generic orbit and the de Rham cohomology of the manifold. To be more precise, our result is valid for a compact Riemannian manifold $M$ on which a connected compact Lie group $G$ acts equidimensionally (i.e. smoothly with all the orbits of the same dimension) and establishes the second page of a spectral sequence, which converges to the de Rham cohomology of $M$ in terms of the former two cohomology theories. In stark contrast to the considerations in \cite{My4}, where the action is assumed to be locally free, in this case already a number of technical difficulties arise due to a potentialy varied cohomological type of the orbits. We provide a number of results which show how our main theorem simplifies under some additional assumptions. Among these results, Theorem \ref{RSSSC} is perhaps the most appealing, due to a vast simplification obtained with the additional assumption that $M$ is simply connected which is very much in the spirit of the often quoted simplified version of the Leray-Serre spectral sequence (with the trivial action of the fundamental group of the base space on the fiber). We also present a method for using our sequence when the action is not equidimensional using an adaptation of the real algebraic blow up to smooth manifolds with a Lie group action. We apply this method to give a topological obstruction to the existence of effective group actions on $4$ and $5$ dimensional manifolds and link it to the classical problem of descerning which manifolds admit a metric with positive scalar curvature.
\newline\indent Let us start by extracting some additional data from our assumptions (and establishing notation in the process).  If we assume that $M$ is connected then the action has generic orbit type $G\slash H$ for some subgroup $H$ which is unique up to conjugation (from now on we reserve $H$ to denote a chosen subgroup in the conjugation class of subgroups corresponding to the generic orbit type). Firstly, let us notice that up to conjugation for each $G$-orbit type $G\slash G_x$ we have that the connected component $H_0$ of the identity in $H$ is also the connected component of identity in $G_x$ (this follows from the Slice Theorem applied to the $G$-action and its equidimensionality). We will denote the corresponding Lie algebras of $G$ and $H$ by  $\mathfrak{g}$ and $\mathfrak{h}$ respectively.
\newline\indent Let us denote by $N$ the normalizer of $H_0$ in $G$ and by $N_0$ the connected component of identity in $N$. Moreover, let $\tilde{N}$ be the subgroup of $N$ preserving a chosen connected component $M^{H_0}_0$ of the space $M^{H_0}$ of fixed points of the action of $H$ on $M$. Then our main result is:
\begin{tw}\label{MainGH} Let $M^{n}$ be a compact connected manifold with an equidimensional action of a compact connected Lie group $G$.  Under the above notation, there is a spectral sequence $E^{p,q}_r$ converging to $H^{\bullet}_{dR}(M)$ with:
$$E^{p,q}_2\cong H^{p}(M_1\slash\mathcal{F}_{1},H^q(\mathfrak{g},\mathfrak{h}))^{N\slash N_0}\cong H^{p}(M_2\slash\mathcal{F}_{2},H^q(\mathfrak{g},\mathfrak{h}))^{\tilde{N}\slash N_0},$$
Where $(M_1,\mathcal{F}_1)$ (resp. $(M_2,\mathcal{F}_2)$) are foliated manifolds covering $(M,\mathcal{F})$ (with $\mathcal{F}$ given by the orbits of the action) with corresponding deck transformation group $N\slash N_0$ (resp. $\tilde{N}\slash N_0$).
\end{tw}
The coverings $M_i$ can be thought of as passage to local coefficients. For the Leray-Serre sequence this passage can be conducted by a universal covering of the base space and taking the cohomology invariant with respect to deck transformations. In our case the orbit space of the $G$ action might not be a manifold and the quotient map to the orbit space might have fibers with different cohomological type. We omit these technical difficulties by choosing to construct appropriate covering spaces for $M$ itself denoted $M_i$. One can think of $M_2$ as the smallest covering spaces that unravels $M$ so that the orbits have cohomology type $H^{\bullet}(\mathfrak{g},\mathfrak{h})$ and the corresponding action of the fundamental group on $H^{\bullet}(\mathfrak{g},\mathfrak{h})$ is trivial. The space $M_1$ consists of multiple copies of $M_2$. The obvious drawback of working with $M_1$ is the lack of connectedness of this spaces. The trade off is the fact that the expression of the second page using $M_1$ somewhat simplifies the groups involved since $N$ and $N_0$ (unlike $\tilde{N}$) depend only on the groups $G$ and $H$.
\newline\indent The paper is organized as follows. In the subsequent section we give some general preliminary material on foliations, Lie groups and Lie algebra cohomology for the readers convienience.
\newline\indent Section $3$ is dedicated to the proof of our main result. In its first subsection we present a special type of Riemannian metric on the manifold $M$ so that the group acts by isometries and when restricted to each $G$ orbit it carries some nice properties similar to that of a metric induced on the quotient by the bi-invariant metric. This enables us to conduct local arguments by freely using the Slice Theorem with the added benefit of the local splitting of the tangent bundle which comes from the global splitting. We follow this in the next subsection by presenting the construction of $M_i$ along with some properties of these covering spaces. The last two subsections are focused on constructing our spectral sequence, establishing its limiting term and use the above covering spaces to compute its first and second page. The computation of the first page consists of three main steps. Firstly, we give a local computation in a neighbourhood of a $G$-orbit. We follow this with a Mayer-Vietoris argument allowing us to compute the first page of the sequence in question for $M_i$. Finally, by passing to parts invariant with respect to the $\tilde{N}\slash N_0$ and $N\slash N_0$ actions on $M_i$ we describe the first page of our spectral sequence as invariant basic forms with coefficients in the Lie algebra cohomology of the pair $(\mathfrak{g},\mathfrak{h})$. We finish the computation by proving that the differential on the first page vanishes on these coefficients. We do this by again first considering the coverings $M_i$ and then passing to invariant parts.
\newline\indent In Section $4$ we present a number of consequences and simplified versions of our spectral sequence under additional assumptions. These can be viewed both as ready made corollaries taylored for applications or as a blue print for the user on how to simplify the sequence for the case of particular interest to them. We split this section into four parts depending on the nature of the simplifications contained in it. In the first part we are focused on topological conditions on $M$ or the space of generic $G$ orbits which simplify the second page by allowing us to omit the coverings. In particular, Theorem \ref{RSSSC} seems to be a good simplification with natural additional assumptions. In the subsequent part we pose a similar question but our focus this time is on the properties of the acting group which might guarantee the desired simplification. We provide here two cases in which $U(n)$ actions (with varying $H_0$) can guatantee the omission of the covering spaces. In the third part we study a Gysin-like long exact sequence for particularly nice choices of $G$ and $H_0$. In the final part we consider the situation of actions with low codimension and arrive at long exact sequences which mimic the Wang exact sequence for bundles over $\mathbb{S}^n$. These form a generalization of similar results from \cite{My4}.
\newline\indent In the final section we present a method of applying our spectral sequence in the non-equidimensional case. We use this method to obtain a result excluding the existence of a non-trivial $\mathbb{S}^3$ action on a large class of $4$-manifolds. While this particular result can be also derived from the classification of $4$-manifolds admitting such actions given in \cite{MP} we extend the methods used in this case to the subsequent subsection in order to give a similar criterion for $5$-manifolds. We also underline the significance of such applications by recalling a well known result stating that if there is no non-trivial action of $\mathbb{S}^3$ on $M$ then there is no effective action of any non-abelian compact connected Lie group on $M$. Hence, our result establishes a reliable way of excluding the existence of such Lie group actions on a given manifold. We finish this section by recalling a link between the existence of such Lie group actions on $M$ and the existence of a metric with positive scalar curvature on $M$ exhibiting potential application to this classical problem.
\newline\indent The paper is written in such a way that a reader mainly interested in applications of our spectral sequence has a decent chance of finding what he needs in the penultimate section of this article. If however, the full power of the spectral sequence is needed it should suffice to read through the explanation of the setup provided in this section together with the main result and possibly parts of Section $3$ describing the construction of $M_i$.

\bigskip
{\Small{\textbf{Acknowledgements} The author would like to thank Professor Martin Saralegui-Aranguren for many detailed conversations on the subject. His insight and remarks proved invaluable during the work on this article.
\newline\indent The work of the author has been partially supported by the Polish National Science Center project number 024/08/X/ST1/01213}}
\bigskip
\section{Preliminaries}
\subsection{Foliations}
We provide a quick review of foliations and basic cohomology.
\begin{defi} A codimension $q$ foliation $\mathcal{F}$ on a smooth $n$-manifold $M$ is given by the following data:
\begin{itemize}
\item An open cover $\mathcal{U}:=\{U_i\}_{i\in I}$ of M.
\item A q-dimensional smooth manifold $T_0$.
\item For each $U_i\in\mathcal{U}$ a submersion $f_i: U_i\rightarrow T_0$ with connected fibers (these fibers are called plaques).
\item For all intersections $U_i\cap U_j\neq\emptyset$ a local diffeomorphism $\gamma_{ij}$ of $T_0$ such that $f_j=\gamma_{ij}\circ f_i$
\end{itemize}
The last condition ensures that plaques glue nicely to form a partition of $M$ consisting of immersed submanifolds of $M$ of codimension $q$. This partition is called a foliation $\mathcal{F}$ of $M$ and the elements of this partition are called leaves of $\mathcal{F}$.
\end{defi}
We call $T=\coprod\limits_{U_i\in\mathcal{U}}f_i(U_i)$ the transverse manifold of $\mathcal{F}$. The local diffeomorphisms $\gamma_{ij}$ generate a pseudogroup $\Gamma$ of transformations on $T$ (called the holonomy pseudogroup). The space of leaves $M\slash\mathcal{F}$ of the foliation $\mathcal{F}$ can be identified with $T\slash\Gamma$.
\begin{defi}
 A smooth form $\omega$ on $M$ is called transverse if for any vector field $X\in\Gamma (T\mathcal{F})$ (where $T\mathcal{F}$ denotes the bundle tangent to the leaves of $\mathcal{F}$) it satisfies $\iota_X\omega=0$. Moreover, if additionally $\iota_Xd\omega=0$ holds for all $X$ tangent to the leaves of $\mathcal{F}$ then $\omega$ is said to be basic. Basic $0$-forms will be called basic functions henceforth.
\end{defi}
Basic forms are in one to one correspondence with $\Gamma$-invariant smooth forms on $T$. It is clear that $d\omega$ is basic for any basic form $\omega$. Hence, the set of basic forms of $\mathcal{F}$ (denoted $\Omega^{\bullet}(M\slash\mathcal{F})$) is a subcomplex of the de Rham complex of $M$. We define the basic cohomology of $\mathcal{F}$ to be the cohomology of this subcomplex and denote it by $H^{\bullet}(M\slash\mathcal{F})$. A transverse structure to $\mathcal{F}$ is a $\Gamma$-invariant structure on $T$. Among such structures the following is most relevant to our work:
\begin{defi}
$\mathcal{F}$ is said to be Riemannian if $T$ has a $\Gamma$-invariant Riemannian metric. This is equivalent to the existence of a Riemannian metric $g$ (called the transverse Riemannian metric) on $N\mathcal{F}:=TM\slash T\mathcal{F}$ with $\mathcal{L}_Xg=0$ for all vector fields $X$ tangent to the leaves.
\end{defi}
We also recall a special class of Riemannian foliations with particularly nice cohomological properties which will be useful in the penultimate section of this article:
\begin{defi} A codimension $q$ foliation $\mathcal{F}$ on a compact connected manifold $M$ is called homologically orientable if $H^{q}(M\slash\mathcal{F})=\mathbb{R}$. A foliation $\mathcal{F}$ on a compact manifold $M$ is called homologically orientable if its restriction to each connected component of $M$ is homologically orientable.
\end{defi}
We will later see that there is a large class of foliations considered in this paper which are homologically orientable. There is also a version of differential operator theory and Hodge theory for manifolds with foliations and basic forms which allows the proof of the folowing result (cf. \cite{E1}):
\begin{tw} Let $\mathcal{F}$ be a Riemannian foliation on a compact manifold $M$. Then:
\begin{enumerate}
\item $H^{\bullet}(M\slash\mathcal{F})$ is finite dimensional.
\item If in addition $\mathcal{F}$ is homologically orientable then there is an isomorphism between $H^{k}(M\slash\mathcal{F})$ and $H^{q-k}(M\slash\mathcal{F})$ (this is the basic version of Poincar\'e duality).
\end{enumerate}
\end{tw}
We are also going to need a result on homotopy invariance of basic cohomology (cf. \cite{BR-A,HR}). We say that a smooth map $f:(M_1,\mathcal{F}_1)\to (M_2,\mathcal{F}_2)$ is foliated if $f$ takes the leaves of $\mathcal{F}_1$ to leaves of $\mathcal{F}_2$. It is easy to see that foliated maps induce a morphism in basic cohomology. Two foliated maps $f_0,f_1:(M_1,\mathcal{F}_1)\to (M_2,\mathcal{F}_2)$ are foliated homotopic if there is a continuous map $H:[0,1]\times M_1\to M_2$ such, that $H(0,x)=f_0(x)$, $H(1,x)=f_1(x)$ and for each $t_0\in [0,1]$ the function $H(t_0,x)$ is foliated (in particular smooth).
\begin{tw}\label{HI}(Theorem $1$ \cite{BR-A}) If $f_0,f_1:(M_1,\mathcal{F}_1)\to (M_2,\mathcal{F}_2)$ are foliated homotopic then the maps induced in basic cohomology by $f_0$ and $f_1$ are equal.
\end{tw}
We finish this section by recalling the spectral sequence of a Riemannian foliation.
\begin{defi}\label{Fili} We put:
$$F^p_{\mathcal{F}}\Omega^{p+q}(M):=\{\alpha\in\Omega^{p+q}(M)\text{ }|\text{ } \iota_{X_{q}}...\iota_{X_0}\alpha=0,\text{ for } X_0,...,X_{q}\in\Gamma(T\mathcal{F})\}.$$
An element of $F^p_{\mathcal{F}}\Omega^{p+q}(M)$ is called an $(p+q)$-differential form of filtration $p$.
\end{defi}
The definition above in fact gives a filtration of the de Rham complex. Hence, via known theory from homological algebra we can construct a spectral sequence as follows:
\begin{enumerate}
\item The $0$-th page is given by $E_0^{p,q}=F_{\mathcal{F}}^{p}\Omega^{p+q}(M)\slash F_{\mathcal{F}}^{p+1}\Omega^{p+q}(M)$ and $d^{p,q}_0:E_0^{p,q}\to E^{p,q+1}_0$ is simply the morphism induced by $d$.
\item The $r$-th page is given inductively by:
$$E_{r}^{p,q}:=Ker(d^{p,q}_{r-1})\slash Im(d^{p-r,q+r-1}_{r-1})=\frac{\{\alpha\in F^{p}_{\mathcal{F}}\Omega^{p+q}(M)\text{ }|\text{ } d\alpha\in F_{\mathcal{F}}^{p+r}\Omega^{p+q+1}(M)\}}{F^{p+1}_{\mathcal{F}}\Omega^{p+q}(M)+d(F_{\mathcal{F}}^{p-r+1}\Omega^{p+q-1}(M))}$$
\item The $r$-th coboundary operator $d_r:E_r^{p,q}\to E^{p+r,q-r+1}_r$ is again just the map induced by $d$ (due to the description of the $r$-th page this has the target specified above and is well defined).
\end{enumerate}
Furthermore, since the filtration is bounded this spectral sequence converges and its final page is isomorphic to the cohomology of the cochain complex (in this case the de Rham cohomology of $M$).
\begin{rem} The above spectral sequence can be thought of as a generalization of the Leray-Serre spectral sequence in de Rham cohomology to arbitrary Riemannian foliations (as opposed to fiber bundles). However, in full generality no good description of the second page exists to the best of our knowledge.
\end{rem}
\subsection{Lie algebras and Lie groups}
In this subsection we give a brief overview of properties of Lie groups and Lie algebras which will be useful in the rest of the article. Throughout this section let $G$ be an $s$-dimensional Lie group and let $\mathfrak{g}$ be the Lie algebra of left invariant vector fields on $G$ which can be classically identified with $T_eG$ (where $e\in G$ denotes the identity element). Choose a basis $\{\tilde{\xi}_1,...,\tilde{\xi}_s\}$ of $\mathfrak{g}$ and let $\{Z_1,...,Z_s\}$ be a set of right invariant vector fields on $G$ such that $(Z_i)_e=(\tilde{\xi}_i)_e$ then it is well known that:
\begin{enumerate}
\item $Z_i$ are the fundamental vector fields of the left action of $G$ on itself,
\item $[\tilde{\xi}_i,\tilde{\xi}_j]_e=-[Z_i,Z_j]_e$ for all $i,j\in \{1,...,s\}$,
\item $[\tilde{\xi}_i,Z_j]=0$ for all $i,j\in \{1,...,s\}$.
\end{enumerate}
Next we are also going to need the following version of the Slice Theorem (cf. \cite{B}):
\begin{tw} Let $G$ be a compact Lie group acting smoothly on a manifold $M$. For each orbit $Gx$ there exists a tubular neighbourhood $U$ equivariantly diffeomorphic to $G\times_{G_x} N_x(Gx)$ such that the image $S$ of elements represented by $\{e\}\times N_x(Gx)$ through this diffeomorphism is a $G_x$-invariant manifold passing through $x$ (where we consider the normal space $N_x(Gx)$ of the orbit of $x$ at $x$ with the action of $G_x$ given by derivatives at $x$). In this case we say that $S$ is a slice through $x$ and that $U$ is a neighbourhood admitting a slice.
\end{tw}
Assuming that $G$ is compact any orbit of a smooth action of $G$ on a manifold $M$ is equivariantly diffeomorphic to $G\slash G_x$ where $G_x$ is the isotropy group at a point $x$ in the given orbit. Moreover, two orbits $Gx$ and $Gy$ are equivariantly diffeomorphic if and only if $G_x$ is conjugate to $G_y$ (in this case we say that the two orbits are of the same type; cf. \cite{B}). Taking $K_i$ to be the union of all orbits of the given type we get a Whitney stratification of $M$ (in the sense of \cite{Str}) with a finite number of strata.
\newline\indent Our main object of study is the following special class of smooth actions:
\begin{defi} A smooth action of $G$ on a smooth manifold $M$ is called equidimensional if all its orbits have the same dimension.
\end{defi}
If the action of $G$ on $M$ is equidimensional then the connected components of orbits form a foliation on $M$. If moreover $G$ is assumed to be connected then this foliation is simply given by the orbits of the action. If additionaly $M$ is connected then by the Slice Theorem there is an open dense subset of $M$ consisting of orbits of the type $G\slash H$ (where $H$ is a fixed representative of the conjugacy class). Furthermore, since the action is equidimensional again by using the Slice Theorem one can see that for an orbit $Gx\cong G\slash G_x$ of an arbitrary point $x$ the conjugacy class of $(G_x)_0$ has to coincide with that of $H_0$ and in fact if we choose the point $x$ such that $H\subset G_x$ then $H_0=(G_x)_0$.
\newline\indent We wish to now recall the construction of the relative Lie algebra cohomology of $(\mathfrak{g},\mathfrak{h})$ and relate it to the basic cohomology of $(G,\mathcal{H})$ (where $\mathcal{H}$ is the foliation given by the orbits of the right action of a connected component of a closed subgroup $H\subset G$ on $G$). The cochain complex, $(\bigwedge\text{}^{\bullet}\mathfrak{g}^*,d_{\mathfrak{g}})$
where,
$$d_{\mathfrak{g}}\alpha(X_0,...,X_{n}):=\sum_{1\leq i<j\leq n}(-1)^{i+j}\alpha([X_i,X_j],X_0,...,\hat{X}_i,...,\hat{X}_j,...,X_n),$$
is called the Chevalley–Eilenberg complex and its cohomology $H^{\bullet}(\mathfrak{g})$ the Lie algebra cohomology of $\mathfrak{g}$. While this cohomology theory has found many applications throughout mathematics its most important property for our purpose is the following string of isomorphisms (for a compact connected Lie group $G$):
$$H^{\bullet}_{dR}(G)\cong H^{\bullet}_{G}(G)\cong H^{\bullet}(\mathfrak{g}),$$
where the middle term is the cohomology of the cochain complex $\Omega_G^{\bullet}(G)$ of left invariant forms on $G$. The first isomorphism is induced simply by the inclusion of forms and can be proven by finding a cochain map that induces the inverse mapping in cohomology (averaging a form over $G$). The second isomorphism follows simply from comparing $d_{\mathfrak{g}}$ with the exterior derivative $d$ restricted to left invariant forms.
A similar thing happens for all homogenous spaces $G\slash H_0$, where $G$ is a compact Lie group and $H_0$ is a closed connected subgroup. For this we need the so called relative Lie algebra cohomology of the pair  $(\mathfrak{g},\mathfrak{h})$. More, precisely consider the following subcomplex of the Chevalley-Eilenberg complex (this can be defined for any pair of Lie algebras $\mathfrak{h}\subset\mathfrak{g}$):
$$\Omega^{\bullet}(\mathfrak{g},\mathfrak{h}):=\{\omega\in\Exterior^{\bullet}\mathfrak{g}^*\text{ }|\text{ }\forall_{X\in\mathfrak{h}}\iota_X\omega=\iota_Xd\omega=0\}.$$
Similarly as before, when $G$ is compact and $H_0$ is closed its cohomology describes the cohomology of the considered homogenous space by the following chain of isomorphisms:
$$H^{\bullet}_{dR}(G\slash H_0)\cong H^{\bullet}_{G}(G\slash H_0)\cong H^{\bullet}(\mathfrak{g},\mathfrak{h}):=H^{\bullet}(\Omega^{\bullet}(\mathfrak{g},\mathfrak{h})).$$
However, if $H_0$ is not assumed to be closed in $G$ then the proper geometric intuition for this cohomology is provided by the basic cohomology of the foliation $\mathcal{H}$ defined by the orbits of the action of $H_0$:
$$H^{\bullet}(G\slash \mathcal{H})\cong H^{\bullet}_{G}(G\slash \mathcal{H})\cong H^{\bullet}(\mathfrak{g},\mathfrak{h}).$$
However, in our context this generalization serves mostly as a link between the basic cohomology on $G$ and the cohomology of the orbit, since all the isotropy groups $G_x$ (and their connected components) have to be closed subgroups of $G$.
\newline\indent In order to describe more general homogenous spaces $G\slash H$ with $H$ possibly not connected (and $H_0$ now denoting the connected component of $H$) an invariant part with respect to an apropriate group of transformations (acting on the right) has to be taken. This group can be of course taken to be $H\slash H_0$ but it can be also considered as a subgroup of $N\slash H_0$ where $N$ denotes the normalizer of $H_0$. This group coincides with the group of equivariant diffeomorphisms of $G\slash H_0$ and hence quotienting out a subgroup of it makes it so that $G$ still acts on the resulting manifold. Furthermore, when considered in cohomology the action of $N\slash H_0$ factors through the action of $N\slash N_0$ (essentially due to the homotopy invariance of basic cohomology presented earlier) and consequently we can describe the cohomology of such homogenous spaces as $H^{\bullet}(\mathfrak{g},\mathfrak{h})^{\hat{N}\slash N_0}$ where $\hat{N}$ is the smallest group containing both $N_0$ and $H$.
\subsection{Classification of orbits for $SU(2)$-actions}
In the final section of this article we analyze $\mathbb{S}^3\cong SU(2)$ actions on $4$-manifolds in order to provide a restricitive obstruction to the existence of effective actions of compact connected non-abelian Lie groups on such manifolds. The argument relies on a partial classification of possible homogeneous spaces (and consequently possible orbit types) for this group. More precisely, we are interested in distinguishing the various homogenous spaces with isotropy groups of dimension strictly greater then $0$. This in turn is mostly straightforward computation using the Lie algebra $\mathfrak{su}(2)$ of $SU(2)$. For the readers convenience we present here this partial classification.
\newline\indent  Let us start by recalling that the Lie algebra of $SU(2)$ is a $3$-dimensional vector space with bracket given by the cross product (since it is isomorphic to the Lie algebra of $SO(3)$).  To classify homogeneous spaces it suffices to classify closed Lie subgroups and this is most conveniently done by first classifying Lie subalgebras of the corresponding Lie algebra. But from the above description of the bracket it is clear that the only non-trivial Lie subalgebras of this Lie algebra are $1$-dimensional subspaces. Consequently, the considered homogenous spaces are either $2$-dimensional or $\mathbb{S}^3\slash\Gamma$ for some finite group $\Gamma$ or trivial.
\newline\indent Let us give some further detail on the $2$-dimensional homogeneous spaces. As can be seen from the above description any connected $1$-dimensional subgroup has to be a circle (otherwise it would have to be a non-closed subgroup isomorphic to $\mathbb{R}$ and consequently its closure would be an abelian subgroup of dimension strictly greater then $1$ providing a contradiction with the previous paragraph). Consequently, these circles form all the $1$-parameter subgroups of $SU(2)$. Furthermore, using the description of $SU(2)$ as unit quaternions $\mathbb{S}^3$ we can see that the quotient of $SU(2)$ by any of these circles is simply $\mathbb{S}^2$ since it corresponds to quotietning out the multiplication by unit complex numbers in $\mathbb{S}^3\subset\mathbb{H}=\mathbb{C}^2$ with the complex structure given by a chosen unit imaginary quaternion. The resulting left action of $\mathbb{S}^1$ on $\mathbb{S}^2$ is simply rotation around a chosen axis twice. To see this note that the central element $-I$ is an element of any of the circles above and consequently the action of $SU(2)$ factors through the double cover map to a transitive action of $SO(3)$ on $\mathbb{S}^2$ with each circle forming a double cover of the corresponding circle in $SO(3)$. The circles in $SO(3)$ are precisely the $1$-parameter subgroups of $SO(3)$ and are well known to act on $\mathbb{S}^2\cong SO(3)\slash\mathbb{S}^1\cong SU(2)\slash\mathbb{S}^1$ by rotating along some axis once. 
\newline\indent If we assume that $H$ is a not necessarily connected, compact, $1$-dimensional subgroup of $\mathbb{S}^3$ then by the above paragraph $\mathbb{S}^3\slash H_0\cong\mathbb{S}^2$. The left action of $H_0$ on $\mathbb{S}^2$ is also specified above and it has two fixed points. On this space $H$ has to act (on the right) by equivariant diffeomorphism for the $\mathbb{S}^3$-action to be preserved. However, such diffeomorphisms have to take $H_0$ fixed points to $H_0$ fixed points and (since the action is transitive) are specified by an image of a single point. Hence, there are only two such diffeomorphism and consequently the action of $H$ factors through the $\mathbb{Z}_2$ action (given by the antipodal map). Moreover, the kernel of this factorization has to be simply $H_0$ since otherwise the $H$-orbit in $\mathbb{S}^3$ and the $H$-orbit in $\mathbb{S}^2\cong\mathbb{S}^3\slash H_0$ would have a different number of connected components which gives a contradiction. Finally, this implies that the only possible two dimensional orbits are $\mathbb{S}^2$ and $\mathbb{R}P^2$.
\newline\indent To end this section let us note that this also implies that the only non-trivial closed normal subgroups of $\mathbb{S}^3$ are finite. To see this it suffices to show that in the classification above neither $\mathbb{S}^2$ nor $\mathbb{R}P^2$ are the underlying smooth manifold of a Lie group. This however is obvious since neither of them is parallelizable (since both have non-zero Euler characteristic). This implies that any non-trivial action of $\mathbb{S}^3$ factors through an effective action of $\mathbb{S}^3\slash\Gamma$ for some finite group $\Gamma$.
\section{Proof of the Main Theorem}
\subsection{A Choice of Riemannian Metric}
In this section, we construct an invariant Riemannian metric appropriate for our purposes. To do this let us fix a bi-invariant metric $g_G$ on $G$ and a cover $\mathcal{U}$ of $M$ by open saturated sets $U_{\alpha}$ admitting a slice $D_{\alpha}$ through a given point $x_{\alpha}$ with isotropy group containing $H$. Let us consider the manifolds $G\times D_{\alpha}$ together with any metric $g_{\alpha}$ such, that for any $(h,x)\in G\times D_{\alpha}$:
\begin{enumerate}
\item $T(\{h\}\times D_{\alpha})$ and $T(G\times\{x\})$ are orthogonal,
\item $g_{\alpha}|_{T(G\times\{x\})}$ is the chosen bi-invariant metric $g_G$.
\item $g_{\alpha}|_{T(\{h\}\times D_{\alpha})}$ is invariant under the action of $G_{x_{\alpha}}$ on the second component. 
\end{enumerate}
Due to the above construction the metric $g_{\alpha}|_{TG_{x_{\alpha}}^{\perp}}$ descends to $U_{\alpha}$ (here by $TG_{x_{\alpha}}$ we mean the bundle tangent to the orbits of the action of $G_{x_{\alpha}}$ simultaneously on the right on $G$ and on the left on $D_{\alpha}$). By taking a $G$-invariant partition of unity $\rho_{\alpha}$ subordinate to $\mathcal{U}$ (such a partition always exists by subjugating any partition of unity subordinate to $\mathcal{U}$ to an averaging process) we can define our new metric $g:=\sum_{\alpha}\rho_{\alpha}g_{\alpha}$. To see the value of this metric let us note that every metric $g_{\alpha}$ restricted to a choosen orbit is precisely the metric induced from the chosen bi-invariant metric on $G$ via the quotient by the isotropy group of a chosen point (different choices of the point will result in the same metric since the initial metric is bi-invariant). This follows since gluing $M$ from $U_{\alpha}$ has to be done by equivariant diffeomorphism (for there to be a group action on $M$) which preserve the given metric on the orbit (since they themselves when restricted to an orbit can be described by the right multiplication by elements from the normalizer of the relevant isotropy group). Since $g$ is a convex combination of such metrics then restricted to the orbits it also has to be induced from that metric (and has to be invariant). Let us state this property:
\begin{prop}\label{Mcomplement} Let $M$ be a compact Riemannian manifold with an equidimensional action of a compact connected Lie group $G$ with a chosen bi-invariant metric $g_G$. Then there exists an invariant Riemannian metric on $M$ such that its restriction to any orbit coincides with a metric induced on that orbit from $g_G$.
\end{prop}
\begin{rem}\label{Perp} In fact we can infer a slightly stronger property which we will use in the computation of the second page. Namely, for any neighbourhood $U$ admitting a perpendicular slice through a given point $x$ on the central orbit we can construct an invariant metric on $G\times D$ which induces the given metric on $U$ such that it agrees on $D$ with the metric induced from the slice, it is the given bi-invariant metric on $G$ and $\{e\}\times D$ is perpendicular to $G\times\{0\}$. To construct such a metric it suffices to lift the given metric on $U$ to a metric on the bundle normal to the orbits of the $G_x$-action and then extend it to $TG_x$ by adding a tensor which on this bundle is the restriction of the given bi-invariant metric and which vanishes when one of the vectors is from the bundle tangent to $\{h\}\times D$ or the orhogonal complement of this bundle with respect to the chosen bi-invariant metric in $G\times\{y\}$. Under our construction the resulting metric on $G$ is the given bi-invariant metric and in fact the constructed metric on $G\times D$ is bi-invariant when restricted to $T(G\times D)|_{Gx}$ (with respect to the action of $G$ on the first component). While it seems like a minor improvement it will prove important in order to consider vector fields perpendicular to an orbit.
\end{rem}
\subsection{A Simplifying Lift for Equidimensional Actions}
Let $(M,g)$ be a compact connected Riemannian manifold on which a compact connected Lie group $G$ acts equidimensionally, where the metric $g$ is as constructed in the previous subsection. In this subsection, we are going to study mostly the topology of $M$ and construct appropriate covering spaces $M_1$ and $M_2$ which will aid us in presenting the second page of our spectral sequence in terms of basic cohomology. The property we are interested in is that the action of the fundamental group $\tilde{\pi}:=\pi_1((M_2\backslash\Sigma_2)\slash G)$ on the generic orbits of the $G$-action (which constitute a fiber bundle) factors through $N_0$ (where $\Sigma_2$ denotes the set of orbits of non-generic type in $M_2$). This will later allow us to describe the second page of the spectral sequences of invariant forms of such a covering as a tensor product.
\newline\indent We begin the construction of our covering spaces by providing an alternate description of the manifold $M$ itself.
\begin{tw} Let $G$ be a connected compact Lie group acting equidimensionally on $M$. Then under the above notation there is an equivariant diffeomorphism:
$$f:G\times_{N}M^{H_0}\to M,$$
where $M^{H_0}$ denotes the points of $M$ which are fixed by the action restricted to $H_0$ and $N$ denotes the normalizer of $H_0$. 
\begin{proof} The diffeomorphism $f$ is given by taking $[g,x]$ and sending it to $gx$. Firstly, let us note that by equidimensionality of the action for each isotropy group $(G_x)_0=H_0$ (up to conjugation). This implies that the action of $N$ on $M^{H_0}$ is well defined. To see this it suffices to prove that this action is well defined on points fixed by $H_0$ in each orbit. This is obvious when the orbit in question is of the type $G\slash H_0$. For more general orbits it suffices to consider the quotient map $p:G\slash H_0\to G\slash G_x$ (with $G_{x0}=H_0$) and note that $p^{-1}([G\slash G_x]^{H_0})=[G\slash H_0]^{H_0}$. From this we infer that the map $f$ is well defined and equivariant.
\newline\indent We shall now prove bijectivity in two steps by first proving that the sets of orbits are send bijectively to each other and then that this map restricted to each orbit is bijective. Firstly, the map between the orbit spaces is surjective since $M^{H_0}$ intersects each orbit of the $G$ action on $M$. Secondly, to prove injectivity it suffices to note that in both spaces the points in $M^{H_0}$ belong to the same $G$ orbit if and only if they differ by an element of $N$. As before it is easy to see for $G\slash H_0$ and can be extended to a general orbit $G\slash G_x$ by the equality $p^{-1}([G\slash G_x]^{H_0})=[G\slash H_0]^{H_0}$. Bijectivity (and in fact equivariant diffeomorphism), of each orbit follows from the fact that for $x\in M^{H_0}$ we get $H_0=(G_x)_0$ and the equality  $p^{-1}([G\slash G_x]^{H_0})=[G\slash H_0]^{H_0}=N\slash H_0$ (where $p$ is as above). Using these facts we get:
$$G\times_{N}(Gx)^{H_0}\cong G\times_{N}(G\slash{G_x})^{H_0}\cong G\times_{N}(G\slash{G_x})^{(G_x)_0}\cong  G\times_{N} N\slash{G_x} \cong G\slash{G_x}\cong Gx.$$
The proof can be now finished by proving that the specified map is in fact a local diffeomorphism by comparing $G\times_N U^{H_0}$ and $U$ (where $U\cong (G\times D)\slash G_x$ is a neighbourhood satisfying the Slice Theorem for the $G$-action on $M$ and $D$ is the disk in the normal bundle). Then the restriction of $f$ factors as follows:
$$G\times_N U^{H_0}\cong G\times_{N}((G\times D)\slash{G_x})^{H_0}\cong G\times_{N}((G\times D)\slash{G_x})^{(G_x)_0}\cong  G\times_{N} (N\times D)\slash{G_x} \cong (G\times D)\slash{G_x}\cong U,$$
where the penultimate diffeomorphism stems from the fact that in $G\times N\times D$ the action of $N$ on the first two terms and the action of $G_x$ on the latter two terms commute.
\end{proof}
\end{tw}
We can further simplify this construction by taking one of the connected components $M^{H_0}_0$ of $M^{H_0}$ which results also in reducing $N$ to a subgroup $\tilde{N}$ of elements which preserve the connected component $M^{H_0}_0$. Explicitly we get:
$$M\cong G\times_{\tilde{N}}M^{H_0}_0.$$
This suggests two sensible candidates for a covering space "unraveling" the $\tilde{N}$ (or $N$) action:
$$M_1\cong G\times_{N_0}M^{H_0},$$
$$M_2\cong G\times_{N_0}M^{H_0}_0.$$
As mentioned in the introduction the covering space $M_2$ exchanges the drawback of introducing an additional group $\tilde{N}$ dependent on the topology of the Lie group action for working with a connected covering space when compared to $M_1$. Let us show that these coverings have indeed the desired property:
\begin{tw}\label{Cover} The map $p_2:M_2\to M_2\slash\tilde{N}\cong M$ is a $G$ equivariant finite covering of $M$ such that the cohomology of each orbit of the action on $M_2$ is isomorphic to $H^{\bullet}(\mathfrak{g},\mathfrak{h})$. Moreover, the action of $\tilde{\pi}$ on the generic orbits of $M_2$ can be represented by elements in $N_0$.
\begin{proof} The map $p_2$ is obviously equivariant. Since $\tilde{N}$ acts freely on $G$ on the right the action of $\tilde{N}\slash{N_0}$ is free on $M_2$ (the action considered here comes from the right action on $G$ and left action on $M^{H_0}$ in the product). Consequently, since $M\cong M_2\slash (\tilde{N}\slash N_0)$ with $\tilde{N}\slash N_0$ a finite group, the map $p_2$ is a finite covering of $M$.
\newline\indent The cohomology of the orbits is a direct computation arising from the study of the given coverings near a given orbit with the use of the Slice Theorem. More explicitly, let $x\in M^{H_0}_0\subset M$ be a point on a given orbit. From the Slice Theorem we know that $Gx$ is a further quotient $G\slash H_1$ of the generic orbit type $G\slash H$ (possibly $H=H_1$) with $H_1\subset N$. Moreover, due to the structure of $M^{H_0}$ the group $H_1$ can be chosen to be a subgroup of $\tilde{N}$. To see this note that it suffices to show that a single point from the component $M^{H_0}_0$ remains in that component under the action of $H_1$. Hence, by choosing $H_1$ (within the appropriate conjugation class inside $N$) so that $M^{H_1}\cap M^{H_0}_0\neq\emptyset$ (such a choice can be made since $M^{H_0}\supset M^{H_1}\neq\emptyset$ and $N$ can transform any component of $M^{H_0}$ to any other component) we get the desired result.
\newline\indent This implies that the cohomology of the orbit through $x$ has to be of the form $H^{\bullet}(\mathfrak{g},\mathfrak{h})^{N_1\slash N_0}$ where $N_1\subset \tilde{N}$ is the group spanned by $N_0$ and $H_1$ (again for this we note that $N_0$ acts trivially on the right on the cohomology of $G\slash H_0$). The covering $p_2$ has degree $|\tilde{N}\slash N_0|$ while the preimage of $Gx$ has $|\tilde{N}\slash N_1|$ connected components each constituting a $|N_1\slash N_0|$-fold covering of $Gx$. Moreover, by applying this construction to the chosen orbit we get that the connected component of $G\times_{N_0} (G\slash H_1)^{H_0}\cong G\times_{N_0}N\slash H_1$ is:
$$G\times_{N_0}N_1\slash H_1\cong G\times_{N_0}N_0\slash (N_0\cap H_1)\cong G\slash (N_0\cap H_1),$$
Hence, the connected covering is just $p:G\slash (N_0\cap H_1)\to G\slash H_1$ with deck transformations given by the action of $H_1\slash (N_0\cap H_1)\cong N_1\slash N_0$ on the right. This implies that the cohomology of the orbits in $M_2$ is of the desired form (since $N_0\cap H_1\subset N_0$ acts trivially on $H^{\bullet}(\mathfrak{g},\mathfrak{h})$). 
\newline\indent For the final statement let us note that it is equivalent to there being no path in $( [\{e\}\times M^{H_0}_0]\backslash\Sigma_2)\subset M_2$ connecting two different connected components of $Gx\cap  [\{e\}\times M^{H_0}_0]$ for some element $x\in [\{e\}\times M^{H_0}_0]\backslash\Sigma_2\subset M_2$. To see this it suffices to trivialize the bundle using neighbourhoods with slices through points in $[e\times M^{H_0}_0]$ and note that the transition between such trivialization are represented by elements of $N$ (or more precisely elements of the noramlizer of $N_0\cap H$ which is a subset of $N$). However, this is obvious due to construction since factoring out the action of $N_0$ in $G\times M^{H_0}_0$ means the intersection $[\{e\}\times M^{H_0}_0]\cap Gx$ has only one connected component in $M_2$.
\end{proof}
\end{tw}
\begin{rem} In the proof we explicitly identify the equivariant diffeomorphism type of orbits in $M_2$. Explicitly, if the image of the considered orbit through $p_2$ is of the form $G\slash H_1$ with $H_1\subset \tilde{N}$ then the covering restricted to that orbit is just $p:G\slash (N_0\cap H_1)\to G\slash H_1$ with deck transformations given by the action of $H_1\slash (N_0\cap H_1)\cong N_1\slash N_0$ on the right (where $N_1\subset\tilde{N}$ is the group generated by $H_1$ and $N_0$).
\end{rem}
By repeating this argument for each connected component of $M_1$ we readily get a similar statement for this space:
\begin{cor}\label{Cover2} The map $p_1:M_1\to M_1\slash N\cong M$ is a $G$ equivariant finite covering of $M$ such that the cohomology of each orbit of the action on $M_1$ is isomorphic to $H^{\bullet}(\mathfrak{g},\mathfrak{h})$. Moreover, the action of $\tilde{\pi}$ on the orbits of generic type of any connected component of $M_1$ can be represented by elements in $N_0$.
\end{cor}
For the sake of brevity, since treating the cases of $M_1$ comes down to repeating the arguments for all the connected components we won't adress them seperately unless more essential changes need to be made in the argument. To use this setting effectively we will need the following observations:  
\begin{obs}\label{Pass}\begin{itemize}
\item Under the above notation $\Omega^{\bullet}_G(M)\cong\Omega^{\bullet}_{G}(M_1)^{N\slash N_0}$ where the action is by the group of deck transformations $N\slash N_0$.
\item Under the above notation $\Omega^{\bullet}_G(M)\cong\Omega^{\bullet}_{G}(M_2)^{\tilde{N}\slash N_0}$ where the action is by the group of deck transformations $\tilde{N}\slash N_0$.
\end{itemize}
\begin{proof} This is immediate from the theory of manifold covering spaces.
\end{proof}
\end{obs}
These observations will allow us to study the cohomology of $M$ using the covering spaces $M_1$ and $M_2$.
\subsection{Construction and the First Page of the Spectral Sequence}
Let us denote by $\mathcal{F}$ the regular foliation given by the orbits of the $G$-action. Since this action is isometric with respect to the chosen metric $g$ these foliations are Riemannian and the bundle $T\mathcal{F}$ and its orthogonal complement bundle $T\mathcal{F}^{\perp}$ are preserved by the action. The splitting:
$$TM=T\mathcal{F}\oplus T\mathcal{F}^{\perp},$$
induces a bigrading in the space of forms $\Omega^{\bullet} (M)$ and consequently on the space of $G$-invariant forms $\Omega^{\bullet}_{G}(M)$. Moreover, since $T\mathcal{F}$ is involutive the operator $d$ when restricted to invariant forms with bidegree $(p,q)$ (denoted $\Omega_G^{p,q} (M)$) splits into three parts $d^{0,1}$, $d^{1,0}$ and $d^{2,-1}$ shifting the degree respectively by $(0,1)$, $(1,0)$ and $(2,-1)$. This follows easilly from the formula:
$$d\alpha(X_0,...,X_k)=\sum (-1)^i\mathcal{L}_{X_i}(\alpha(X_0,...,\hat{X}_i,...,X_k))+\sum_{i<j} (-1)^{i+j}\alpha([X_i,X_j],X_0,....,\hat{X_i},...,\hat{X_j},...,X_k).$$
\newline\indent The spectral sequence we are considering arises from the complex $\Omega^{\bullet}_{G}(M)$ endowed with the restriction of the filtration from Definition \ref{Fili} i.e.:
$$F^p_{\mathcal{F}}\Omega_{G}^{p+q}(M):=\{\alpha\in\Omega_{G}^{p+q}(M)\text{ }|\text{ } \iota_{X_{q}}...\iota_{X_0}\alpha=0,\text{ for } X_0,...,X_{q}\in\Gamma(T\mathcal{F})\}.$$
The cohomology of this complex (and hence also the limit of the spectral sequence) is isomorphic with de Rham cohomology of $M$. Hence, to prove the main Theorem it suffices to compute the second page of this spectral sequence.
\newline\indent In order to deal with the above spectral sequence efficiently it will be useful to relate the above filtration to the bigrading introduced earlier. For this we identify $\Omega^{p,q}(M)$ with $E^{p,q}_0$ via the obvious isomorphism of vector spaces (each equivalence class in $E_0^{p,q}$ has a unique representative in $\Omega^{p,q}_G(M)$). Moreover, to compute the second page let us note that similarly as for a spectral sequence corresponding to a double complex the operator $d_0^{p,q}$ is just $d^{1,0}$ on $\Omega^{p,q}_G(M)$ and $d_1^{p,q}$ is just the operator induced in $E_1^{p,q}$ by $d^{0,1}$ on $\Omega^{p,q}_G(M)$. Again analogously as in the case of a double complex (in fact even more so since $d^{2,-1}$ might not be trivial in our case) the description of $d_r^{p,q}$ in terms of this splitting become somewhat more complicated for $r\geq 2$ but for our purpose this is sufficient.
\newline\indent Passing to the covering spaces let us note that due to the local character of the splitting and the filtration if we repeat the same construction for $M_1$ and $M_2$ then the maps induced by $p_1$ and $p_2$ preserve both the splitting and the filtration. Hence, instead of studying $\Omega_{G}(M)$ we can study $\Omega_{G}(M_1)$ and $\Omega_{G}(M_2)$ (with their filtration and splitting) so that passing to the cohomology of $\Omega_{G}(M)$ can be done by taking the part invariant with respect to $\tilde{N}\slash N_0$ or $N\slash N_0$. Throughout the proof we are going to stick to the notation $\mathcal{F}_{i}$ for the foliation on $M_i$ by orbits of the action.
\newline\indent To compute the first page of the spectral sequence for $M_2$ let us first pick a saturated neigbourhood $V$ of an arbitrary point $x\in M_2$ which admits a slice $D$ through $x$ perpendicular to $Gx$ and compute the first page of the spectral sequence for the complex $\Omega^{\bullet}_{G}(V)$ (considered with an analogue of the above filtration).
\begin{lem}\label{Loc} With notation as above the first page $E_1^{p,q}(V)$ of the above spectral sequence is isomorphic to $\Omega^{p}(V\slash (\mathcal{F}_2|_V),H^q(\mathfrak{g},\mathfrak{h}))$.
\begin{proof} Firstly, let us lift the problem to $G\times D$ from $V\cong (G\times D)\slash G_x$. Let us denote by $\mathcal{H}$ the foliations on $G\times D$ by the connected components of the orbits of the actions of $G_x$. With this notation one can observe that the forms that descend to $V$ (as $G$-invariant forms) are precisely the basic forms on  $(G\times D,\mathcal{H})$ invariant under the actions of $G$ and $G_x$ (where the former acts from the left on $G$ and the later acts simultaneously on the right on $G$ and on the left on $D$). Moreover, note that due to the definition of $\mathcal{H}$ any basic form is already invariant under the action of the identity component $G_{x0}=H_0$ of $G_x$ and hence it is enough to consider $G_x\slash G_{x0}$ in the previous sentence. Hence, we have an isomorphism:
$$\Omega^{\bullet}_{G}(V)\cong[\Omega^{\bullet}_{G}(G\times D\slash\mathcal{H})]^{G_x\slash G_{x0}}.$$
The left action of $G$ on $G\times D$ allows us to lift the foliation $\mathcal{F}_2|_V$ to a foliation $\tilde{\mathcal{F}}_2$ on $G\times D$ given by $G$-orbits (or equivalently with leaves given as the preimages of leaves of $\mathcal{F}_2|_V$). This allows us to pullback the metric from $TV$ to $N\mathcal{H}$ to get a $G$-invariant splitting:
$$N\mathcal{H}\cong (T\tilde{\mathcal{F}}_2\slash T\mathcal{H})\oplus (T\tilde{\mathcal{F}}_2\slash T\mathcal{H})^{\perp} \cong (T\tilde{\mathcal{F}}_2\slash T\mathcal{H})\oplus T(G\times D)\slash T\tilde{\mathcal{F}}_2,$$
From this we see that:
$$\Omega^{\bullet}_{G}(G\times D\slash\mathcal{H})\cong\Omega^{\bullet}((G\times D)\slash (\tilde{\mathcal{F}}_2))\otimes \Omega^{\bullet}(\mathfrak{g},\mathfrak{h})\cong\Omega^{\bullet}((G\times D)\slash (\tilde{\mathcal{F}}_2), \Omega^{\bullet}(\mathfrak{g},\mathfrak{h})),$$
Now by noting that the quotient map preserves the splittings of $N\mathcal{H}$ and $TV$ we can compute the first page of the desired spectral sequence from the first page of the spectral sequence corresponding to the obvious filtration of the complex $\Omega^{\bullet}((G\times D)\slash \tilde{\mathcal{F}}_2, \Omega^{\bullet}(\mathfrak{g},\mathfrak{h}))$. Hence, since $d^{0,1}$ vanishes on basic forms (since on basic forms $d=d^{1,0}$) whereas restricted to invariant forms on orbits its just the operator $d_{\mathfrak{g}}$ and taking cohomology commutes with the action of a finite group on a cochain complex we conclude:
$$E^{p,q}_1(V)\cong H^{p}([\Omega^{\bullet}((G\times D)\slash (\tilde{\mathcal{F}}_2), \Omega^{q}(\mathfrak{g},\mathfrak{h}))]^{G_x\slash G_{x0}},d^{0,1})$$
$$\cong [\Omega^{p}((G\times D)\slash (\tilde{\mathcal{F}}_2), H^{q}(\mathfrak{g},\mathfrak{h}))]^{G_x\slash G_{x0}}.$$
Hence, by Theorem \ref{Cover} we get that $ H^{q}(\mathfrak{g},\mathfrak{h})$ is invariant under the action of $G_x\slash G_{x0}$ and consequently:
$$E^{p,q}_1(V)\cong [\Omega^{p}((G\times D)\slash (\tilde{\mathcal{F}}_2), H^{q}(\mathfrak{g},\mathfrak{h}))]^{G_x\slash G_{x0}}\cong\Omega^{p}(V\slash (\mathcal{F}_2|_V),H^q(\mathfrak{g},\mathfrak{h})).$$
\end{proof}
\end{lem}
Moreover, let us note that a similar result will hold by a similar computation for any saturated open subset of $V$. Finally, we note that since the action of the group $G$ and the metric $g$ are globally defined the transition functions between two saturated open sets as in the previous lemma when restricted to each $G$-orbit have to be given by equivariant isometries. Moreover, by the last statement in Theorem \ref{Cover} and the fact that excising $\Sigma_2$ does not disconnected any neighbourhood admitting a slice in $M_2$, we can cover $M_2$ by open sets such that these isometries can be all represented by elements of $N_0$ on generic orbits. To see that excising $\Sigma_2$ does not disconnect neighbourhoods admitting slices note that this could happen only if there is some strata of codimension $1$. However, for a generic point $x$ in this component $G_x$ has to fix a hyperplane in $N_x(Gx)$ corresponding to this component of $\Sigma_2$. Consequently it can (and must for the orbit to not be of generic type) act non-trivially only on its orthogonal complement. Since, this complement is $1$-dimensional it follows that this action factors through the only possible $O(1)\cong\mathbb{Z}_2$ action making the trace of this strata a $1$-sided submanifold in any neighbourhood admitting a slice around a generic point. Since excising a finite familly of $1$-sided manifolds will not disconnect the neighbourhood we indeed see that excising $\Sigma_2$ cannot do so (by first excising all the strata of lower dimension in the trace of $\Sigma_2$ and then excising this familly).
\newline\indent But since $N_0$ is connected the isometries of the orbits have to be represented by elements from $N_0$ for non-generic orbits as well (by approaching the transition function on a non-generic orbit with ones on nearby generic orbits). Since such maps act trivially on the cohomology of orbits (by a homotopy argument for $H^{\bullet}_{dR}(G\slash H_0)\cong H^{\bullet}(\mathfrak{g},\mathfrak{h})$), we can easilly glue such subsets via a version of the Mayer-Vietoris Theorem for $E^{p,q}_0=\{\Omega_{G}^{p,q}(V),d^{0,1}\}$:
\begin{lem}\label{MV0} Let $(M,g)$ be a compact connected Riemannian manifold with an equidimensional action of $G$. Moreover, let $U$ and $V$ be open saturated subsets. Then the following short sequence of cochain complexes is exact:
$$0\to \{\Omega_G^{p,\bullet}(U\cup V),d^{0,1}\}\to\{\Omega_G^{p,\bullet}(U),d^{0,1}\}\oplus\{\Omega_G^{p,\bullet}(V),d^{0,1}\}\to\{\Omega_G^{p,\bullet}(U\cap V),d^{0,1}\}\to 0.$$
\begin{proof} The only non-trival part is the surjectivity of the last non-trivial chain map. This is proved via the same argument as in the proof of the Mayer-Vietoris Theorem in \cite{ttt} after noting that we can construct a partition of unity subordinate to $\{U,V\}$ consisting of basic functions (take any partition of unity suboordinate to $\{U,V\}$ and average it over $G$).
\end{proof}
\end{lem}
Using the above Lemma along with isomorphisms from Lemma \ref{Loc} (for open saturated subsets of neighbourhoods admitting slices) and the triviality of the transition functions allows us to compare $E^{p,q}_1$ with the basic forms with values in $H^{\bullet}(\mathfrak{g},\mathfrak{h})$ giving us the following ladder diagram:
$${\tiny\begin{tikzcd}
... \arrow[r]& E^{p,q-1}_1(U)\oplus E^{p,q-1}_1(V)\arrow[r]\arrow[d]& E^{p,q-1}_1(U\cap V)\arrow[r]\arrow[d]& E^{p,q}_{1}(U\cup V)\arrow[r]& ...\\
... \arrow[r]& \Omega^{p}(U\slash\mathcal{F}_2,H^{q-1}(\mathfrak{g},\mathfrak{h}))\oplus \Omega^{p}(V\slash\mathcal{F}_2,H^{q-1}(\mathfrak{g},\mathfrak{h}))\arrow[r]& \Omega^{p}(U\cap V\slash\mathcal{F}_2, H^{q-1}(\mathfrak{g},\mathfrak{h}))\arrow[r,"0"]& \Omega^{p}(U\cup V\slash\mathcal{F}_2, H^q(\mathfrak{g},\mathfrak{h}))\arrow[r]& ...
\end{tikzcd}}$$
where $U$ and $V$ are neighbourhoods admitting slices in $M_2$. To define the missing vertical arrow let us first note that from the square on the left and the fact that the lower horizontal arrow is epimorphic the same has to be true for the upper horizontal arrow (since the two vertical arrows are isomorphisms). Which implies that the horizontal arrow above the zero mapping is zero as well. Hence, $E^{p,q}_1(U\cup V)$ and  $\Omega^{p}(U\cup V\slash\mathcal{F}_2, H^{q-1}(\mathfrak{g},\mathfrak{h}))$ are precisely the kernels of the subsequent arrows and hence the map between them is given by the restriction of the subsequent isomorphism $ E^{p,q}_1(U)\oplus E^{p,q}_1(V)\cong \Omega^{p}(U\slash\mathcal{F}_2,H^{q}(\mathfrak{g},\mathfrak{h}))\oplus \Omega^{p}(V\slash\mathcal{F}_2,H^{q}(\mathfrak{g},\mathfrak{h}))$. Consequently, this arrow is an isomorphism as well. By consecutively applying this proceedure to glue one additional subset from the chosen cover at a time we get the first page of our spectral sequence:
$$E_1^{p,q}(M_2)=\Omega^{p}(M_2\slash\mathcal{F}_{2},H^q(\mathfrak{g},\mathfrak{h})).$$
Since the operation of taking an invariant part with respect to a finite group commutes with taking cohomology we get in the end:
$$E_1^{p,q}(M)=\Omega^{p}(M_2\slash\mathcal{F}_{2},H^q(\mathfrak{g},\mathfrak{h}))^{\tilde{N}\slash N_0}.$$
Similar statement follows for $M_1$.
\subsection{The second page of the spectral sequence}
To compute the second page $E^{p,q}_2(M)$ we will first show that locally $d^{1,0}$ vanishes on elements from $\Omega^{\bullet}(\mathfrak{g},\mathfrak{h})^{N_0}$, treated locally as a subspace of invariant forms via the splitting of $TM_2$. More precisely, let $V$ be a neighbourhood of any point $x\in M^{H_0}_2$ admitting a slice $D$ through $x$ perpendicular to $Gx$ and consider the quotient map $\pi:G\times D\to G\times D\slash G_x\cong V$. We consider an invariant Riemannian metric $g$ such as in the previous subsections. We can naturally identify $\Omega^{\bullet}(\mathfrak{g},\mathfrak{h})^{N_0}$ with a subspace of $\Omega_{G}^{\bullet}(G\times D\slash\mathcal{H})\subset \Omega_{G}^{\bullet}(G\times D)$ using the orthogonal complement of $T\tilde{\mathcal{F}_2}\slash T\mathcal{H}$ in $N\mathcal{H}$ (where $\tilde{\mathcal{F}_2}$, $\mathcal{H}$ and the metric on $N\mathcal{H}$ are as in the previous section). Note that all such forms descend to $\Omega^{\bullet}_G(V)$ since $G_x\subset N_0$ and hence $\Omega^{\bullet}(\mathfrak{g},\mathfrak{h})^{N_0}$ can be identified with a subspace of $\Omega^{\bullet}_G(V)$. Moreover, note that this space does not depend on the choice of slice through $x$ since if we repeat the construction for a different choice it would simply change each orbit by a right multiplication by some element of $N_0$ (possibly different for each orbit). Since, regardless of the choice of slice the descended forms vanish on the orthogonal complement of $T\mathcal{F}_2|_V$ and restricted to the orbit are $N_0$ invariant this shows that indeed this space does not depend on the chosen slice through $x$. Our goal is to show that $d^{1,0}$ vanishes on these descended forms.
\newline\indent We start by working on the dense and open set $U$ consisting of the points belonging to orbits of generic type (equivariantly diffeomorphic to $G\slash (H\cap N_0)$). This set by the Slice Theorem has the structure of a locally trivial bundle over $U\slash G$ with fiber $G\slash (H\cap N_0)$. Let $V$ be a saturated neighbourhood of a point $x\in U$ admitting a slice $D$ through $x$ perpendicular to $Gx$.
\newline\indent We consider a frame $X_1,...,X_s,Y_1,...,Y_q$ for $G\times D$ where:
\begin{itemize}
\item $X_i$ are left invariant vector fields spaning $T\tilde{\mathcal{F}}_2$ which are constant in the direction of $D$.
\item $Y_j$ are left invariant vector fields spaning $(T\tilde{\mathcal{F}}_2)^{\perp}$.
\end{itemize} 
Let us use this frame to compute $d^{1,0}\tilde{\alpha}$ for $\tilde{\alpha}\in\Omega^{\bullet}(\mathfrak{g},\mathfrak{h})^{N_0}$ (the lift of $\alpha$ to $G\times D$). For this we write for $Y\in\{Y_1,...,Y_q\}$:
$$d\tilde{\alpha} (Y,X_{i_1},...,X_{i_k})=\mathcal{L}_{Y}(\tilde{\alpha}(X_{i_1},...,X_{i_k}))+\sum_l (-1)^l\mathcal{L}_{X_{i_l}}(\tilde{\alpha}(Y, X_{i_1},...,\hat{X}_{i_l},...,X_{i_k}))$$
$$+\sum_l (-1)^l\tilde{\alpha}([Y,X_{i_l}],X_{i_1},...,\hat{X}_{i_l},...,X_{i_k})+\sum_{l<r} (-1)^{l+r}\tilde{\alpha}([X_{i_l},X_{i_r}], Y,...,\hat{X}_{i_l},...,\hat{X}_{i_r},...,X_{i_k}),$$
and consider the terms on the right at a lift $\tilde{x}$ of $x$. The first term vanishes since by our choice of $X_i$ the function arising from evaluating $\tilde{\alpha}$ on these vector fields is constant. The second and fourth term vanish since they include the evaluation of $\tilde{\alpha}$ on a perpendicular vector field $Y$. Similarly, the third term has to be zero as well since each $[Y,X_{i_l}]$ is perpendicular (since by Remark \ref{Perp} the lifted metric is bi-invariant on $T(G\times D)|_{G\tilde{x}}$ with respect to the left and right multiplications in $G$). Consequently, $d^{1,0}\tilde{\alpha}=0$ at $\tilde{x}$ which implies that $d^{1,0}\alpha=0$ at $x$ (since the map induced on forms by the quotient $G\times D\to G\slash (H\cap N_0)\times D\cong V$ is pointwise monomorphic) and consequently $d^{1,0}\alpha=0$ at any point in $Gx$. To show that this is the case for other orbits through $y\in V$ we repeat this computation for a slice which is perpendicular to the chosen orbit (we use the fact that the space $\Omega^{\bullet}(\mathfrak{g},\mathfrak{h})^{N_0}$ does not depend on the choice of slice) which can be done since a connected component of the restriction of this slice to an appropriate saturated tubular neighbourhood of the chosen orbit is a slice in its own right. Hence, we have proved that $d^{1,0}$ vanishes locally on elements from $\Omega^{\bullet}(\mathfrak{g},\mathfrak{h})^{N_0}$ at any point $x\in U$. To pass to general orbits it suffices to consider the Slice Theorem for a non-generic orbit and approach that orbit by points in $U$.


This finishes the proof of Theorem \ref{MainGH} since it implies that $d^{1,0}$ on the first page doesn't affect $H^{\bullet}(\mathfrak{g},\mathfrak{h})$. To see this note that for $M_2$ classes in $H^{\bullet}(\mathfrak{g},\mathfrak{h})$ can be locally (in the forms on a neighbourhood admitting a slice) represented by elements in $\Omega^{\bullet}(\mathfrak{g},\mathfrak{h})^{N_0}$. Since the transition functions between neighbourhoods admitting slices preserve orthogonal complements and restricted to each $G$ orbit act by an element of $N_0$ these forms can be glued to global representatives of clasess in $\Omega^{0}(M_2\slash\mathcal{F}_2,H^{\bullet}(\mathfrak{g},\mathfrak{h})^{N_0})\cong \Omega^{0}(M_2\slash\mathcal{F}_2,H^{\bullet}(\mathfrak{g},\mathfrak{h}))$ on which $d^{1,0}$ vanishes. Consequently, passing to cohomology (with respect to $d^{1,0}$) will leave the coefficients (identified with the classes represented by the above constructed basic $0$-forms) intact and hence the second page is just:

$$E_2^{p,q}(M)\cong H^{p}(M_2\slash\mathcal{F}_{2},H^q(\mathfrak{g},\mathfrak{h}))^{\tilde{N}\slash N_0}.$$

\begin{rem} Alternatively, one could finish the argument that $d^{1,0}$ doesn't affect the coefficients by noting that the above local pictures glue to a trivial vector bundle $M_2\times H^{\bullet}(\mathfrak{g},\mathfrak{h})$ (since all the transition functions come from $N_0$ and hence are trivial). If we now consider the operator $d^{1,0}$ on basic forms with values in this bundle we can see that it vanishes on constant (with respect to the presented transition function) sections due to the above consideration (since it also can be glued from the local descriptions above).
\end{rem}
It is also interesting that using the description with $M_1$ we can transform the space:
$$M_1\slash\mathcal{F}_1\cong G\times_{N_0}M^{H_0}\slash G\cong M^{H_0}\slash N_0,$$
giving us:
\begin{cor} Let $M^{n+s}$ be a compact connected manifold with an equidimensional action of a compact connected $s$-dimensional Lie group $G$. Under the above notation, there is a spectral sequence $E^{p,q}_r$ converging to $H^{\bullet}_{dR}(M)$ with:
$$E^{p,q}_2\cong H^{p}_{dR}(M^{H_0}\slash N_0,H^q(\mathfrak{g},\mathfrak{h}))^{N\slash N_0}.$$
\end{cor}
\section{Special Cases, Examples and Applications}
\subsection{Under what topological conditions can the coverings be ommited?}
In this section we study some topological conditions under which the main Theorem can be simplified eliminating the need to take into account the action of $\tilde{N}$. These resutls can be thought of as the analogue of the case when the fundamental group of the base acts trivially on the fiber in the Leray-Serre spectral sequence.
\begin{tw} Under the notation and assumptions of Theorem \ref{MainGH} let us assume additionally that the fundamental group of $(M\backslash\Sigma)\slash G$ acts trivially on the cohomology of the generic orbit which is isomorphic to $H^q(\mathfrak{g},\mathfrak{h})$. Then the second page of the spectral sequence of the action is given by:
$$E^{p,q}_2\cong H^{p}(M\slash\mathcal{F},H^q(\mathfrak{g},\mathfrak{h})).$$
\begin{proof} 

From the above asumption we know that the deck transformation group $\tilde{N}\slash N_0$ has to act trivially on the cohomology of the generic orbit in $M_2\backslash\Sigma_2$ (again here we consider the trivial bundle $M_2\times H^{\bullet}(\mathfrak{g},\mathfrak{h})$ like at the end of the previous section). By looking at the first page of the spectral sequence (noting that this action extends to $M_2$) we see that this action has to be trivial on all the orbits by the density of $M_2\backslash\Sigma_2$. Consequently, the action of $\tilde{N}$ on the second page can be factored as follows:

$$H^{p}(M_2\slash\mathcal{F}_{2},H^q(\mathfrak{g},\mathfrak{h}))^{\tilde{N}\slash N_0}\cong H^{p}(M_2\slash\mathcal{F}_{2})^{\tilde{N}\slash N_0}\otimes H^q(\mathfrak{g},\mathfrak{h})\cong H^{p}(M\slash\mathcal{F},H^q(\mathfrak{g},\mathfrak{h})), $$
finishing the proof.
\end{proof}
\end{tw}
While the above considerations mimic the ones in the classical Leray-Serre spectral sequence quite faithfully in our setting they seem somewhat convoluted and artificial. This is mainly due to the fact that for fiber bundles the base space is often the initial point of the consideration whereas in our case the initial point is $M$ and hence instead of placing conditions on some subset of the orbit space it seems to be more natural to place conditions on $M$ itself. The subsequent result shows that this line of thought also gives relevant simplifications of the spectral sequence.
\begin{tw}\label{RSSSC} Under the notation and assumptions of Theorem \ref{MainGH} let us assume additionally that $M$ is simply connected. Then the second page of the spectral sequence is given by:
$$E^{p,q}_2\cong H^{p}(M\slash\mathcal{F},H^q(\mathfrak{g},\mathfrak{h})).$$
\begin{proof} The result is practically immediate by noting that $M_2$ is a connected covering space of $M$ and hence by our additional assumption it has to be equal to $M$.
\end{proof}
\end{tw}
\subsection{For what groups can passing to the coverings be ommitted?}
In the previous subsection we gave topological conditions on the action of $G$ on $M$ which guarantee that the technical tool of passing to covering spaces can be ommited. In this subsection we will study certain special cases in which we can infer the same conclusion based only on the structure of the acting groups $G,H,N,N_0$ (i.e. the set of data that encode the structure of the generic orbit without reference to the group $\tilde{N}$ encoding the global twist and possible singular orbits). Firstly, let us note that if $N$ acts trivially on $H^{\bullet}(\mathfrak{g},\mathfrak{h})$ then the spectral sequence has second page:
$$E^{p,q}_2\cong H^{p}(M\slash\mathcal{F},H^q(\mathfrak{g},\mathfrak{h})).$$
This follows immediately from Theorem \ref{MainGH} since the action on the coefficients is trivial by the above assumpiton while $H^{\bullet}(M_2\slash\mathcal{F}_2)^{\tilde{N}\slash N_0}\cong H^{\bullet}(M\slash\mathcal{F})$. This happens for example when $N$ is connected.
\newline\indent A particular example of this situation is when $G=U(n)$ and $H=U(k)\times \{I_{n-k}\}$ (which we denote simply $U(k)$ in the following consideration). In this case the normalizer of $H$ is $U(k)\times U(n-k)$. To see this let us first note that $U(k)$ is in fact normal in $U(k)\times U(n-k)$. To show that the normalizer cannot be bigger let us take $h\in U(n)\backslash (U(k)\times U(n-k))$ and assume that $h$ normalizes $U(k)$. By the assumption $h$ does not preserve $\mathbb{C}^k\times\{0\}$. Hence, we can choose a vector $v\in\mathbb{C}^n\backslash (\mathbb{C}^k\times\{0\})$ such that $h^{-1}v\in \mathbb{C}^k\times\{0\}$ now choose an element $u\in U(k)$ such that $h^{-1}v$ generates a one dimensional eigenspace of the eigenvalue $\lambda\neq 1$ of $u$. Then $huh^{-1}v=\lambda v$ and consequently $huh^{-1}$ cannot be an element of $U(k)$ which gives the desired contradiction.
\begin{cor} Let $(M,g)$ be a compact connected manifold with an equidimensional action of $U(n)$ with the connected component of the isotropy group of a generic orbit given by $U(k)\times\{I_{n-k}\}$ (up to conjugation). Then, there is a spectral sequence $E^{p,q}_r$ converging to $H^{\bullet}_{dR}(M)$ with:
$$E^{p,q}_2\cong H^{p}(M\slash\mathcal{F},H^q(\mathfrak{u}(n),\mathfrak{u}(k))).$$
\end{cor}
\indent On the other hand, $U(k)\times U(n-k)$ by a very similar argument as in the previous paragraph is self normalizing in $U(n)$ unless $n=2k$ and hence gives us a second example.
\begin{cor} Let $(M,g)$ be a compact connected manifold with an equidimensional action of $U(n)$ with the connected component of the isotropy group of a generic orbit given by $U(k)\times U(n-k)$ (up to conjugation) with $2k\neq n$. Then, there is a spectral sequence $E^{p,q}_r$ converging to $H^{\bullet}_{dR}(M)$ with:
$$E^{p,q}_2\cong H^{p}(M\slash\mathcal{F},H^q(\mathfrak{u}(n),\mathfrak{u}(k)\times\mathfrak{u}(n-k))).$$
\end{cor}
\indent In the excluded case $2k=n$ the normalizer of $U(k)\times U(n-k)=U(k)\times U(k)$ is the group generated by $U(k)\times U(k)$ and switching the places of the first $k$ coordinates with the last $k$ coordinates and so there are two possibilities depending on whether $\tilde{N}=N_0$ or $\tilde{N}=N$ (where in the former case the Theorems remain the same as presented before while in the later case a $\mathbb{Z}_2$ action and a double covering have to be taken into consideration).
\subsection{Gysin long exact sequences} In this section we provide an application of our spectral sequence for a particularly nice choice of $(\mathfrak{g},\mathfrak{h})$. More precisely, if for a given $l\in\mathbb{N}$ we have:
$$H^{i}(\mathfrak{g},\mathfrak{h})=\left\{\begin{array}{l}
\mathbb{R},\text{ for }i\in\{0,l\},\\
0,\text{ otherwise}.
\end{array} \right.$$
We can try to construct a Gysin type long exact sequence. While in full generality the invariant part makes it hard to consider such a sequence, under some conditions which will guarantee that the action of $\tilde{N}$ on the top cohomology is trivial it is a straightforward consequence of our main Theorem. We present two instances of such long exact sequence. In the first case we give the precise conclusion of the above consideration:
\begin{tw} Let $M$ be a compact, connected manifold with an equidimensional action of a compact connected Lie group $G$ with the connected component of the isotropy group of a generic orbit given by $H_0$ (up to conjugation). Moreover, assume that for a given $l\in\mathbb{N}$ we have:
$$H^{i}(\mathfrak{g},\mathfrak{h})=\left\{\begin{array}{l}
\mathbb{R},\text{ for }i\in\{0,l\},\\
0,\text{ otherwise},
\end{array} \right.$$
and that $\tilde{N}$ acts trivially on $H^{i}(\mathfrak{g},\mathfrak{h})$. Then we have a long exact sequence:
$$...\to H^{k}_{dR}(M)\to H^{k-l}(M\slash\mathcal{F})\to H^{k+1}(M\slash\mathcal{F})\to  H^{k+1}_{dR}(M)\to ...$$.
\begin{proof} This is a straightforward consequence of the second page having only two nonzero rows which is guaranteed by the assumptions on the cohomology of the pair $(\mathfrak{g},\mathfrak{h})$. Note that similarly as before the action on the coefficients can be ommited by the above assumption and $H^{\bullet}(M_2\slash\mathcal{F}_2)^{\tilde{N}\slash N_0}\cong H^{\bullet}(M\slash\mathcal{F})$ allowing us to ommit the coverings in the statement.
\end{proof}
\end{tw}
The second is a corollary of the above result paired with the considerations in the first subsection of this section which allows us to replace the technical assumption on $\tilde{N}$ with the assumption that $M$ is simply connected:
\begin{cor} Let $M$ be a compact, connected, simply connected manifold with an equidimensional action of a compact connected Lie group $G$ with the connected component of the isotropy group of a generic orbit given by $H_0$ (up to conjugation). Moreover, assume that for a given $l\in\mathbb{N}$ we have:
$$H^{i}(\mathfrak{g},\mathfrak{h})=\left\{\begin{array}{l}
\mathbb{R},\text{ for }i\in\{0,l\},\\
0,\text{ otherwise}.
\end{array} \right.$$
Then we have a long exact sequence:
$$...\to H^{k}_{dR}(M)\to H^{k-l}(M\slash\mathcal{F})\to H^{k+1}(M\slash\mathcal{F})\to  H^{k+1}_{dR}(M)\to ...$$.
\end{cor}
The above cases are taylored so that the result does not depend on the coverings. Generalizations using covering spaces are of course possible as well. As an example let us go in greater detail to the case $G=SO(l+1)$ and $H_0=SO(l)$, for $l\geq 2$. In this case, the normalizer is $O(l)\subset SO(l+1)$ where the inclusion is given by taking a martix $A\in O(l)$ to a block matrix $diag(A, det(A))$. Hence, the generic orbit is $SO(l+1)\slash SO(l)\cong\mathbb{S}^l$ or $SO(l+1)\slash O(l)\cong\mathbb{R}P^l$. Consequently, $N\slash N_0\cong \mathbb{Z}_2$ where the action of this group on $H^{\bullet}(\mathfrak{g},\mathfrak{h})$ is trivial in all degrees except possibly for the top non-vanishing cohomology on which it acts by multiplication by $(-1)^{l+1}$. Using the description of our spectral sequence for $M_1$ (which is a double cover in this case) and the fact that it has at most two non-vanishing rows we get a long exact sequence:
\begin{tw} Let $M$ be a compact connected manifold with an equidimensional action of $SO(l+1)$ with the connected component of the isotropy group of a generic orbit $SO(l)$. Then we have a long exact sequence:
$$...\to H^{k}_{dR}(M)\to H^{k-l}(M_1\slash\mathcal{F}_1)^{\mathbb{Z}_2}\to H^{k+1}(M\slash\mathcal{F})\to  H^{k+1}_{dR}(M)\to ...,$$
where the action on the basic cohomology above is just  $(-1)^{l+1}$ times the action induced by the deck transformations on basic cohomology. By substituting $M_1$ and noting that $H_0=N_0$ in this case we can present our sequence in an alternate form. Hence, we get a sequence:
$$...\to H^{k}_{dR}(M)\to H^{k-l}_{dR}(M^{H_0})^{\mathbb{Z}_2}\to H^{k+1}(M\slash \mathcal{F})\to  H^{k+1}_{dR}(M)\to ...$$
\end{tw}
\begin{rem} We wish to note that for $l=2$ this sequence is a special case of the main result of \cite{RPSA}.
\end{rem}
Moreover, if $l$ is even then the above sequence has to split:
\begin{cor} Let $M$ be a compact connected manifold with an equidimensional action of $SO(l+1)$ with the connected component of the isotropy group of a generic orbit $SO(l)$ with $l$ even. Then for any $k$ we have short exact sequences:
$$0\to H^{k}(M\slash\mathcal{F})\to H^{k}_{dR}(M)\to H^{k-l}(M_1\slash\mathcal{F}_1)^{\mathbb{Z}_2}\to 0,$$
$$0\to H^{k}(M\slash \mathcal{F})\to H^{k}_{dR}(M)\to H^{k-l}_{dR}(M^{H_0})^{\mathbb{Z}_2}\to 0.$$
\begin{proof}The second familly of short exact sequences is a consequence of the first one as in the above theorem. To prove the exactness of the first family it suffices to prove that $d_{l+1}$ vanishes. To see this we note that on the level of forms ($E_0^{p,q}$) the first row for the sequence on $M_1$ is already zero. This follows since otherwise we would have a non-vanishing invariant form in the first row on $\mathbb{S}^l\times \mathbb{D}^n$ (the neighbourhood of a generic orbit in $M^1$). But this would imply the existence of a non-vanishing invariant form on $\mathbb{S}^l$ (by evaluating the form in $\mathbb{S}^l\times \mathbb{D}^n$ on appropriate invariant vector fields spanning the orthogonal complement of the vertical bundle) and hence a nowhere vanishing $1$-form on $\mathbb{S}^l$ which gives a contradiction. However, since $d=d^{0,1}+d^{1,0}+d^{2,-1}$ this implies $d_{l+1}\alpha$ is represented by an exact basic form (since $d^{1,0}$ is the only non-trivial part of $d$ targeting the $0$-th row). This finishes the proof.
\end{proof}
\end{cor}
\subsection{Wang long exact sequences}
In this section we provide applications of our spectral sequence for low codimension actions. These results for the most part can be thought of as a generalization of similar results from \cite{My4} which in turn generalized the Wang long exact sequence. We start with the following:
\begin{prop}\label{Hor}  Let $M^{n+s}$ be a connected, compact, orientable, manifold with an equidimensional action of a compact connected $s$-dimensional Lie group $G$. Then the foliation induced by the action is homologically orientable.
\begin{proof} It is well known (cf. \cite{E1}) that the top basic cohomology of a Riemannian foliation on a compact connected manifold is either $0$ or $\mathbb{R}$. In this case it cannot be $0$ since then we could compute from our spectral sequence that $H^{n+s}_{dR}(M)\cong E^{n,s}_2=0$ which is a contradiction with the orientability of $M$. Hence, the top basic cohomology is isomorphic to $\mathbb{R}$ which means that the foliation is homologically orientable.
\end{proof}
\end{prop}
\begin{tw} Let $M^{1+s}$ be a compact connected oriented manifold with an equidimensional action of a compac Lie group $G$ with orbits of dimension $s$. Then, $M$ cannot be simply connected.
\begin{proof} Assume to the contrary that $M$ is simply connected. Since $E^{0,q}_2$ and $E^{1,q}_2$ are the only non-zero columns on the second page we conclude that the sequence degenerates at this page. Due to Proposition \ref{Hor} we can use the basic version of Poincar\'e Duality and the fact that $H^0(\mathfrak{g},\mathfrak{h})$ is invariant under the action of the normalizer (since its represented by constant functions on $G\slash H_0$) to conclude from $E_2^{p,q}\cong E^{p,q}_{\infty}$ that:
$$H^{1}_{dR}(M)\cong E^{0,1}_2\oplus E^{1,0}_2\cong E^{0,1}_2\oplus H^{0}(\mathfrak{g},\mathfrak{h}).$$
Moreover, since $H^{0}(\mathfrak{g},\mathfrak{h})\neq 0$ together with the above splitting implies that $H^{1}_{dR}(M)\neq 0$ we get the desired contradiction. \end{proof}
\end{tw}
\begin{tw} Let $M^{2+s}$ be a compact, connected, simply connected and oriented manifold with an equidimensional action of a compact Lie group $G$ with orbits of dimension $s$. Then we have a long exact sequence:
$$...\to H^{k}_{dR}(M)\to H^{k}(\mathfrak{g},\mathfrak{h})\to H^{k-1}(\mathfrak{g},\mathfrak{h})\to  H^{k+1}_{dR}(M)\to ...$$.
\begin{proof} We use the well known fact that the morphism $i^*:H^1(M\slash\mathcal{F})\to H^{1}_{dR}(M)$ induced by the inclusion of forms is a monomorphism to conclude that $H^1(M\slash\mathcal{F})=0$ and hence only the $0$-th and the $2$-nd column of $E_2^{p,q}$ are nonzero. Moreover, due to Proposition \ref{Hor} we can use the basic version of Poincar\'e duality together with Theorem \ref{RSSSC} to conclude that  $E^{0,q}_2\cong H^{q}(\mathfrak{g},\mathfrak{h})$ and $E^{2,q}_2\cong H^{q}(\mathfrak{g},\mathfrak{h})$. Hence, we can use the derivative on the second page to arrive at the above long exact sequence.
\end{proof}
\end{tw}

\begin{tw} Let $M^{3+s}$ be a compact, connected, simply connected and oriented manifold with an equidimensional action of a compact Lie group $G$ with orbits of dimension $s$. Then we have a long exact sequence:
$$...\to H^{k}_{dR}(M)\to H^{k}(\mathfrak{g},\mathfrak{h})\to H^{k-2}(\mathfrak{g},\mathfrak{h})\to  H^{k+1}_{dR}(M)\to ...$$
Moreover, $0\cong H^{1}(\mathfrak{g},\mathfrak{h})$.
\begin{proof} We again use the well known fact that the morphism $i^*:H^1(M\slash\mathcal{F})\to H^{1}_{dR}(M)$ induced by the inclusion of forms is a monomorphism to conclude that $H^1(M\slash\mathcal{F})=0$. Moreover, we can use Proposition \ref{Hor} to conclude from the basic version of Poincar\'e duality that $H^2(M\slash\mathcal{F})=0$ and hence, that  $E^{0,q}_2\cong H^{q}(\mathfrak{g},\mathfrak{h})$ and $E^{3,q}_2\cong H^{q}(\mathfrak{g},\mathfrak{h})$ are the only non-zero columns. Hence, we can use the derivative on the third page to arrive at the above long exact sequence. The final assertion follows from the fact that there are no non-zero coboundary operators with source or target on $E^{0,1}_r$ for $r\geq 2$ and consequently $H^1(\mathfrak{g},\mathfrak{h})\cong E^{0,1}_2\cong E^{0,1}_{\infty}\cong H^1_{dR}(M)=0$.
\end{proof}
\end{tw}
\section{Excluding actions of $\mathbb{S}^3$}
\subsection{A blow-up technique.} In this section, we show how to use a version of (real) blow up in order to produce an equidimensional action of the given Lie group on the blown up manifold. This allows us to mitigate the main difficulty of bypassing the equidimensionality assumption when applying our spectral sequence.
\newline\indent Let $M$ be a compact smooth manifold with a compact smooth submanifold $S$ of codimension $k\geq 2$. Then the real blow up $\hat{M}$ of $M$ along $S$ can be described by excising a saturated tubular neighbourhood $T$ (identified with the open unit disk bundle $D(NS)$ in $NS$) of $S$ from $M$ and performing one of the two operations (yielding diffeomorphic results) below:
\begin{enumerate} 
\item Quotienting the boundary of the resulting manifold with boundary, by a $\mathbb{Z}_2$-action given by the multiplication by $\{1,-1\}$ (inherited from $NS$).
\item Gluing back the fiber product $P_{NS}\times_{O(k)}\hat{\mathbb{D}}^k$ using the obvious diffeomorphism $P_{NS}\times_{O(k)}(\hat{\mathbb{D}}^k\backslash\mathbb{R}P^{k-1})\cong P_{NS}\times_{O(k)}(\mathbb{D}^k\backslash\{0\})\cong \overline{D(NS)}\backslash S$ (allowing us to glue back $P_{NS}\times_{O(k)}(\hat{\mathbb{D}}^k\backslash\mathbb{R}P^{k-1})$ in place of $D(NS)$). Here $P_{NS}$ is the orthonormal frame bundle of $NS$ and $\hat{\mathbb{D}}^k$ denotes the standard real algebraic blow up of the closed $k$-disc $\mathbb{D}^k$ at the origin (which is diffeomorphic to the closed unit disc bundle of the canonical line bundle).
\end{enumerate}
To see that the above descriptions are equivalent we note that the obvious diffeomorphism $NS\backslash S\cong NS\backslash \overline{T}$ (given by contracting the radius) can be extended over the projective spaces to constitute a diffeomorphism between one description to the other (since it preserves the lines). 
\newline\indent Moreover, the result does not depend on the choice of tubular neighbourhood since different choices of tubular neighbourhoods differ by an (ambient) isotopy and a bundle isomorphism. Consequently, since such a bundle isomorphism on a complement of the zero section can be always extended to $P_{NS}\times_{O(k)}\hat{\mathbb{D}}^k$ (since it preserves line), gluing back $P_{NS}\times_{O(k)}\hat{\mathbb{D}}^k$ with respect to two different tubular neighbourhoods will result in diffeomorphic manifolds.
\begin{rem} The above construction is an adaptation to the smooth setting of the version of blowing up found in real algebraic geometry similarly as the complex blow up is adapted for the purposes of e.g. symplectic topology.
\end{rem}
We will now show that if there is a Lie group $G$ acting on $M$ and $S$ is saturated (i.e. $x\in S$ implies $Gx\subset S$) then $\hat{M}$ carries an action of $G$ as well (which extends the action on $M\backslash \overline{T}$). For this let us equip $M$ with an invariant Riemannian metric $g$. Since, $S$ is saturated we have $GS=S$ and consequently the action of $G$ has to preserve the bundle tangent to $S$. However, since the metric is invariant the action has to preserve the bundle perpendicular to $TS$ as well. This gives us an action of $G$ on $TS^{\perp}\cong NS$ which is equivariantly diffeomorphic to a saturated tubular neigbourhood of $S$ by a diffeomorphism taking the image of the zero section to $S$ (given by the exponential map). Moreover, since the $G$ action preserves the chosen Riemannian metric it has to preserve the open unit disc bundle $D(NS)$ in $NS$ and the unit sphere bundle constituting the boundary of its closure. Consequently, we have an action of $G$ on $NS\backslash D(NS)$. Finally, since the action is linear it will survive quotienting out the $\mathbb{Z}_2$-action on the boundary and hence we indeed have a $G$-action on $\hat{M}$.
\newline\indent Let us now show how to use this process in order to make the action equidimensional. For this let us consider the stratification of $M$ given by orbit types and let us start with the strata $K_1$ corresponding to a chosen group $H_1$ such that there is no group $\tilde{H}_1\supset H_1$ such that $G\slash\tilde{H}_1$ is an orbit type which occurs in $M$.
\begin{rem} Note that since $M$ is assumed to be compact there are only finitely many orbit types occuring in it. This follows by covering $M$ with neighbourhoods of each orbit which admit a slice and using compactness to choose a finite subcover. Hence, we can always choose $K_1$ as prescribed above. Moreover, this implies that $K_1$ is a compact submanifold of $M$.
\end{rem}
Let us now blow up $K_1$. Firstly, let us note that excising a tubular neighbourhood of $K_1$ will get rid of the chosen orbit type $G\slash H_1$ but unfortunately in principal new orbit types might arise from the $\mathbb{Z}_2$ quotient. Hence, let us choose an orbit $Gx$ through $x\in\partial D(NK_1)$. If $-x\notin Gx$ then the quotient by the $\mathbb{Z}_2$-action simply identifies two different orbits of the same type and consequently no new orbit types arise from this in $\hat{M}$. Unfortunately, it is easy to see that the other case is both plausible and will inevitably produce orbits of a different type. However, one can note that the dimension of the new orbit has to be the same as that of $Gx$. This leads us to the conclusion that in order to make the action equidimensional one should blow up the set $S_1$ of all the orbits of the smallest dimension instead of just $K_1$. However, such a set is only an immersed submanifold (it might have self intersections) and so the presented blow up technique cannot be applied to it immediately.
\newline\indent Fortunately, this set can be made into a submanifold by consequitvely blowing up closures of strata which are compact manifolds. Let us start by commenting on the fact that the set of orbits of the lowest dimension is indeed a compact immersed submanifold of $M$. This is an immediate consequence of the Slice Theorem. Using it we can see that firstly if an orbit of lowest dimension has no orbits of different type of the same dimension arbitrarily close then near that point $S_1$ is indeed a manifold (this is readilly visible since vectors in $N_x(Gx)$ invariant under $G_x$ form a vector subspace). Secondly, by considering the invariant subspaces in $N_x(Gx)$ corresponding to each isotropy group occuring near the chosen orbit we can conclude that a neighbourhood of such an orbit in $M$ is diffeomorphic to an open set in which the trace of $S_1$ is just a number of smooth manifolds intersecting each other.
\newline\indent Let us now choose $K_1\subset S_1$ as before and study the result blowing it up (denoted $\hat{M}_1$) in $M$ has on $S_1$. For this let us denote by $K_2$ the strata in $M\backslash \overline{T}_1$ (where $T_1$ is the chosen tubular neighbourhood for the blow up) corresponding to a chosen group $H_2$ such that there is no group $\tilde{H}_2\supset H_2$ such that $G\slash\tilde{H}_2$ is an orbit type which occurs in $M\backslash K_1$. By considering, $M\backslash \overline{T}_1$ as a subset in $\hat{M}$ (using the first description) we can take its closure $\overline{K_2}$ in that space (this consists of the corresponding strata in $\hat{M}$ and the possible orbits which are double covered by $G\slash H_2$ which arise on the projective spaces). Our claim is that this space is now a smooth submanifold. Indeed by analyzing the effect blowing up has in neighbourhoods admitting slices we can see that locally (near points in $K_1$) the effect it has on each invariant vector space is just blowing up the vector space corresponding to $K_1$ within that vector space. In particular, $\overline{K}_2$ is an immersed smooth manifold. Moreover, note that this part of the argument works for all the closures of strata from $S_1\subset M$ in $\hat{M}_1$. To see that it is indeed a manifold it suffices to note that there are no points of self intersection. However, this is obvious since such a point would have to lie on an orbit which is a further quotient of $G\slash H_2$. Hence, by our construction and the choice of $H_2$ it would have to lie on the exceptional divisor (i.e. the image in $\hat{M}_1$ from the boundary in $M\backslash T_1$). However, by looking at the neighbourhoods admitting slices in $M$ again, taking into account that the strata correspond to invariant subspaces and the fact that by our choice of $H_2$ the only intersection points of the vector spaces corresponding to these orbits lies in the excised parts, we can see that there can be no intersection points on the exceptional divisor as well. Continuing, this process inductively blowing up the closures of the strata corresponding to orbit types occuring in $M$ we will be eventually placed (for each connected component of $S_1$) in the situation where there remains only one orbit type from $M$ of the given dimension with all the remaining orbit types of that dimension (in this component) lying in its closure which is a smooth manifold. Hence, blowing up the closures of these connected components we get rid of the orbits of the lowest dimension. Then we repeat the process for all the dimensions strictly lower then the dimension of the generic orbit.
\newline\indent Our hope now is that in certain cases this process coupled with passing to covers will change the cohomology of the space in a fashion predictable enough so that from applying our sequence to this blown up manifold we can succesfully study the properties of manifolds with arbitrary actions of a given compact group. We give an interesting example of this in the subsequent section along with some further consequences.


\subsection{Excluding actions of $\mathbb{S}^3$ on $4$-manifolds.} The purpose of this subsection is to prove a cohomological criterion excluding the existence of $\mathbb{S}^3$-actions on $4$-manifolds (without relying on the classification from \cite{MP}) by applying our spectral sequence together with the blow up construction to showcase how to effectively combine these methods.
\begin{tw} Let $M$ be a compact connected smooth $4$-manifold with $dim(H^2_{dR}(M))\geq 3$. Then there is no non-trivial (smooth) action of $\mathbb{S}^3\cong SU(2)$ on $M$.
\end{tw}
By a simple corollary (which we state explicitly in the final subsection as Proposition \ref{NAA2}) of Proposition 3.1 from Chapter IV of \cite{Spin} (which we recall in the final subsection as Proposition \ref{NAA}) this implies the following more general result.
\begin{cor} Let $M$ be a compact connected smooth $4$-manifold with $dim(H^2_{dR}(M))\geq 3$. Then there is no effective (smooth) action of any compact connected non-abelian Lie group $G$ on $M$.
\end{cor}
To prove the above theorem let us analyze possibilities depending on the generic orbit type. Note that the only orbit types for the action of $\mathbb{S}^3$ are $\mathbb{S}^3\slash\Gamma$ for some finite subgroup $\Gamma$, $\mathbb{S}^2$, $\mathbb{R}P^2$ and fixed points.
\newline\indent When the generic orbit type is $\mathbb{S}^3\slash\Gamma$ we get two cases:
\begin{itemize}
\item There are no different orbit types. In this case $M$ is a bundle with fiber $\mathbb{S}^3\slash\Gamma$ over a circle. This satisfies the assumption on cohomology by the Leray-Serre spectral sequence and the fact that the de Rham cohomology of $\mathbb{S}^3\slash\Gamma$ is the same as that of $\mathbb{S}^3$ (which follows from computing the fundamental group of $\mathbb{S}^3\slash\Gamma$ and using Poincar\'{e} Duality). For this note that in this case the action of the fundamental group on the fiber has to be trivial since for there to be an action on the resulting bundle the transition has to be done by equivariant diffeomorphisms which are given by right multiplication by some elements of $\mathbb{S}^3\cong SU(2)$ (and these act trivially on $H^{\bullet}(\mathfrak{su}(2))\cong H^{\bullet}(\mathfrak{su}(2))^{\Gamma}\cong H^{\bullet}_{dR}(\mathbb{S}^3\slash\Gamma)$).
\item There are orbits of different type in $M$. Hence, by connectedness the set consisting of points on generic orbits is a disjoint union of copies of $(\mathbb{S}^3\slash\Gamma)\times (0,1)$.
\end{itemize}
Consequently, it suffices to study possible orbits which might "fit" at the ends of the interval in the second case. Note that indeed as we move along the interval to $0$ or $1$ we necessary have to approach a single orbit (since to each value of the interval there is only a single orbit assigned). Hence, we need to study the various possibilities arising from the slice Theorem.
\begin{enumerate}
\item Assuming that the orbit in the middle is $3$-dimensional its neighbourhood admitting a slice is given by $\mathbb{S}^3\times \mathbb{R}\slash \hat{\Gamma}$. Since $O(1)\cong\mathbb{Z}_2$ and by our assumption on the generic orbit the only possible scenario is that $\hat{\Gamma}\slash\Gamma\cong\mathbb{Z}_2$ acts on $\mathbb{R}$ by $-1$. This follows since if the action of the group $G_x$ on $\mathbb{R}$ was trivial then the middle orbit would be generic while an action of a larger group would have a different generic orbit type. Note that this singular orbit cannot join two copies of $(\mathbb{S}^3\slash\Gamma)\times (0,1)$.
\item Assuming that the orbit in the middle is $2$-dimensional there are two possibilities. Firstly, it can be diffeomorphic to $\mathbb{S}^2$ when $\Gamma$ is a finite cyclic subgroup of a circle. The second case (orbit diffeomorphic to $\mathbb{R}P^2$) is similar with $\Gamma$ a finite subgroup of the group $G_x$ consisting of two circles. Neither of these can join two copies of $(\mathbb{S}^3\slash\Gamma)\times (0,1)$.
\item Assuming that the orbit in the middle is $0$ dimensional leaves us with only one possibility for a neighbourhood admitting a slice (a $4$-disc with an obvious action of $\mathbb{S}^3$). Again note that this singular orbit cannot join two copies of $(\mathbb{S}^3\slash\Gamma)\times (0,1)$.
\end{enumerate}
Hence, by the connectedness assumption and the above analysis we can conclude that the set of points on generic orbits is just one copy of $(\mathbb{S}^3\slash\Gamma)\times (0,1)$. We can now finish the case when the generic orbit is $\mathbb{S}^3\slash\Gamma$ by showing via a Mayer-Vietoris argument that no combination of two of the possible orbit types above will make it so that $dim(H^2_{dR}(M))\geq 3$. This is apparent since the cohomology of the intersection is just that of $\mathbb{S}^3\slash\Gamma$ (which again by Poincar\'{e} Duality is the same as that of $\mathbb{S}^3$) and so its first and second cohomology vanish. This implies that $H^2_{dR}(M)$ is just the sum of the second cohomology of the two neighbourhoods admitting slices (where each of them is at most $1$-dimensional by the above classification).
\newline\indent Before jumping in to the most technically involved case let us quickly note that the case when the generic orbit is $0$-dimensional obviously contains only the trivial action and so there is nothing to prove here.
\newline\indent For the case when the generic orbit is $2$-dimensional we will use a combination of the exclusion argument from the $3$-dimensional case and resolution of non-generic orbits applied to the remaining possible orbit types. More precisely, in this case the only possible non-generic orbits are copies of $\mathbb{R}P^2$ (if the generic orbit is $\mathbb{S}^2$) and fixed points. However, one can note that since what is required is the equidimensionality of the action the $\mathbb{R}P^2$ orbits do not have to be resolved. To resolve the strata of fixed points we will need the following slightly technical lemma:
\begin{lem} Let $M$ be a compact connected smooth $4$-manifold with a (smooth) $\mathbb{S}^3$-action with generic orbit of dimension $2$. Then the fixed point strata consists of a disjoint union of circles.
\begin{proof} Firstly, let us note that using the fact that for fixed points there are no orbit types which are further quotients the fixed point strata has to form a compact submanifold. Moreover, since the generic orbit has dimension $2$ the rank of the normal bundle of the fixed point strata has to be at least $3$. Hence, the strata of fixed points has to have dimension either $0$ or $1$.
\newline\indent Let us now assume to the contrary that $x_0$ is an isolated fixed point. Moreover, let us choose an $\mathbb{S}^3$-invariant metric on $M$. For this metric the group action induces an action on an arbitrarily small sphere around $x_0$. Since, $x_0$ is isolated, on sufficiently small spheres there can be no fixed points. Hence, the induced action on such spheres is equidimensional (all the orbits have dimension $2$). However, this leads to a contradiction since by our spectral sequence such an action on a sphere cannot exist. More precisely, since the sphere $\mathbb{S}^3$ is connected and simply connected the second page of our spectral sequence in this case by Theorem \ref{RSSSC} is just:
$$E^{p,q}_2\cong H^{p}(\mathbb{S}^3\slash\mathcal{F}, H^{q}(\mathbb{S}^2)).$$
However, since only two subsequent columns are non-zero on the second page it has to be the last page as well. This gives the desired contradiction since it would imply that the second cohomology of $\mathbb{S}^3$ is non-zero.
\end{proof}
\end{lem}
Hence, we can conclude from the lemma that in order to make the action (with generic orbit of dimension $2$) equidimensional it suffices to blow up circles consisting of fixed points. Let us now show that this process preserves de Rham cohomology. Ideologically this is due to the fact that $\mathbb{S}^1$ and $\mathbb{S}^1\times\mathbb{R}P^2$ have the same cohomology with real coefficients (by the K\"unneth Formula).
\begin{lem} Let $M$ be a compact connected smooth $4$-manifold with a (smooth) $\mathbb{S}^3$-action with generic orbit of dimension $2$ which contains a circle consisting of fixed points. Let us denote by $\hat{M}$ the blow up of that circle in $M$. Then $H^{\bullet}_{dR}(M)\cong H^{\bullet}_{dR}(\hat{M})$.
\begin{proof} Let us consider the smooth map $f:\hat{M}\to M$ (blow down) defined using the second description of the blow up process as follows. On $V:=M\backslash \overline{T}$ (where $T$ is the tubular saturated neighbourhood of the circle) it is simply the identity. On $\hat{U}:=P_{NS}\times_{O(3)}\hat{\mathbb{R}}^3\to P_{NS}\times_{O(3)}\mathbb{R}^3=:U$ it is given by identity on the first coordinate and collapsing $\mathbb{R}P^2$ to the origin in $\mathbb{R}^3$ on the second coordinate. For this map let us write down the Mayer-Vietoris ladder diagram:
$$\begin{tikzcd}
... \arrow[r]& H^{k-1}_{dR}(U)\oplus H^{k-1}_{dR}(V)\arrow[r]\arrow[d]& H^{k-1}_{dR}(U\cap V)\arrow[r]\arrow[d]& H^{k}_{dR}(M)\arrow[r]\arrow[d]&H^{k}_{dR}(U)\oplus H^{k}_{dR}(V)\arrow[r]\arrow[d]& ...\\
... \arrow[r]& H^{k-1}_{dR}(\hat{U})\oplus H^{k-1}_{dR}(V)\arrow[r]& H^{k-1}_{dR}(\hat{U}\cap V)\arrow[r]& H^{k}_{dR}(\hat{M})\arrow[r]&H^{k}_{dR}(\hat{U})\oplus H^{k}_{dR}(V)\arrow[r]& ...
\end{tikzcd}$$
We are going to use the five lemma to show that $f$ induces an isomorphism in cohomology. For this let us note that the map on $V$ is the identity and hence induces an isomorphism. Moreover, $f$ restricts to a diffeomorphism from $\hat{U}\cap V$ to $U\cap V$. Hence, it suffices to show that $f|_{\hat{U}}:\hat{U}\to U$ induces an isomorphism in cohomology. By homotopy invariance it suffices to prove that this is true for $f|_{P_{NS}\times_{O(3)}\mathbb{R}P^{2}}:P_{NS}\times_{O(3)}\mathbb{R}P^{2}\to\mathbb{S}^1$. For this let us note that $P_{NS}\times_{O(3)}\mathbb{R}P^{2}$ is just $\mathbb{S}^1\times\mathbb{R}P^2$. This is visible since it has to be an $\mathbb{R}P^2$ bundle over the circle which is glued by equivariant diffeomorphisms of $\mathbb{R}P^2$. However, the only equivariant diffeomorphism of $\mathbb{R}P^2\cong \mathbb{S}^3\slash H$ (with respect to the $\mathbb{S}^3$-action) is the identity. To see this one can note that $H$ is just a copy of two circles in $\mathbb{S}^3$ since $\mathbb{S}^3\slash\mathbb{S}^1\cong\mathbb{S}^2$. This group acting on $\mathbb{R}P^2$ has a unique fixed point (corresponding to the axis of the chosen rotation on $\mathbb{S}^2$ spanning the circle). Consequently, any such diffeomorphism has to preserve this fixed point. However, since the $\mathbb{S}^3$-action on the orbit is transitive an equivariant diffeomorphism is determined by the image of a single point (and consequently there can be only one such diffeomorphism). Hence, $f|_{P_{NS}\times_{O(3)}\mathbb{R}P^{2}}$ is just the projection in the trivial $\mathbb{R}P^2$ bundle over $\mathbb{S}^1$ and consequently it induces an isomorphism in cohomology finishing the proof by the five lemma.
\end{proof}
\end{lem}
Let us finish the proof of the Theorem by considering the covering space $\hat{M}_2$ of $\hat{M}$ (which is the result of blowing up all the circles of fixed points in $M$). Since it is a finite covering of $\hat{M}$ it induces a monomorphism in cohomology by the description of the cohomology of $\hat{M}$ as invariant (with respect to the deck transformations) cohomology classes on $\hat{M}_2$. Consequently, since $\hat{M}$ and $M$ have the same second cohomology it suffices to show that $dim(H^2_{dR}(\hat{M}_2))\leq 2$.
\newline\indent To see that this is the case note that for a four dimensional manifold with an action of $\mathbb{S}^3$ with only $2$-dimensional orbits the non-zero entries on the second page can only be situated in a $3$ by $3$ square ($p,q\in\{0,1,2\}$). Moreover, one notes that the middle row is neccesarily zero since in this case $H^1(\mathfrak{g},\mathfrak{h})\cong H^{1}_{dR}(\mathbb{S}^2)=0$. Consequently, the sequence for $\hat{M}_2$ degenerates at the second page and hence the dimension of the second cohomology of $\hat{M}_2$ is the sum of the dimensions of $E^{0,2}_2$ and $E^{2,0}_2$. The former of these is just $H^2(\mathfrak{g},\mathfrak{h})\cong\mathbb{R}$ since $\hat{M}_2$ is connected. For the later note that since $\mathbb{S}^3$ is compact the foliation by orbits is Riemannian and hence its top basic cohomology $H^{2}(M\slash\mathcal{F})\cong E_2^{2,0}$ has to be either $\mathbb{R}$ or $0$. Hence, indeed if $\mathbb{S}^3$ acts on a $4$-manifold then the dimension of its second cohomology cannot exceed $2$ finishing the proof.
\subsection{Excluding actions of $\mathbb{S}^3$ on $5$-manifolds.} In this subsection we take our methods a step further to give a similar (although slightly more technical) criterion for excluding actions of $\mathbb{S}^3$ on $5$-manifolds.
\begin{tw}\label{5M} Let $M$ be a compact smooth $5$-manifold. Then:
\begin{enumerate}
\item If $M$ admits an effective (smooth) action of $\mathbb{S}^3$ with generic orbit of dimension $3$ then there is a hyperplane in $H_2(M,\mathbb{R})$ generated by spheres.
\item If $M$ admits an effective (smooth) action of $\mathbb{S}^3$ with generic orbit of dimension $2$ then there is a hyperplane in $H^2_{dR}(M)$ restricted to which the cup product vanishes.
\end{enumerate}
\end{tw}
As in the previous section this can be generalized to arbitrary non-abelian Lie group action in the following way using Proposition \ref{NAA2}:
\begin{cor} Let $M$ be a compact smooth $5$-manifold admitting an effective (smooth) action of a compact connected non-abelian Lie group. Then at least one of the following conditions has to hold:
\begin{enumerate}
\item There is a hyperplane in $H_2(M;\mathbb{R})$ generated by spheres.
\item There is a hyperplane in $H^2_{dR}(M)$ restricted to which the cup product vanishes.
\end{enumerate}
\end{cor}
The remainder of this subsection is dedicated to the proof of Theorem \ref{5M}. This time let us firstly consider the case when the generic orbit is of dimension $2$ due to its similarities to the corresponding argument in the previous subsection.
\begin{lem} Let $M$ be a compact connected smooth $5$-manifold with a (smooth) $\mathbb{S}^3$-action with generic orbit of dimension $2$. Then the fixed point strata consists of a disjoint union of compact surfaces.
\begin{proof} Firstly, let us note that using the fact that for fixed points there are no orbit types which are further quotients the fixed point strata has to form a compact submanifold. Moreover, since the generic orbit has dimension $2$ the rank of the normal bundle of the fixed point strata has to be at least $3$. Hence, the strata of fixed points has to have dimension either $0$, $1$ or $2$.
\newline\indent Let us now assume to the contrary that $x_0$ is an isolated fixed point. Moreover, let us choose an $\mathbb{S}^3$-invariant metric on $M$. For this metric the group action induces an action on an arbitrarily small sphere around $x_0$. Since, $x_0$ is isolated, on sufficiently small spheres there can be no fixed points. Hence, the induced action on such spheres is equidimensional (all the orbits have dimension $2$). However, this leads to a contradiction since by our spectral sequence such an action on a sphere cannot exist. More precisely, since the sphere $\mathbb{S}^4$ is connected and simply connected the second page of our spectral sequence in this case by Theorem \ref{RSSSC} is just:
$$E^{p,q}_2\cong H^{p}(\mathbb{S}^4\slash\mathcal{F}, H^{q}(\mathbb{S}^2)).$$
However, since the second page consists of only two non-zero rows with non-zero terms only for $p\in\{0,1,2\}$ seperated by a row of zeroes it has to be the last page as well. This gives the desired contradiction since it would imply that the second cohomology of $\mathbb{S}^4$ is non-zero.
\newline\indent To exclude the existence of a circle consisting of fixed points we can give a similar argument. In this case for the chosen metric consider a geodesic $4$-disc passing through the chosen point $x_0$ in the given circle of fixed points which is perpendicular to that circle. The boundary sphere $\mathbb{S}^3$ has to be preserved by the action due to invariance of the metric and hence $\mathbb{S}^3$ would have to admit an action of the group $\mathbb{S}^3$ with all orbits of dimension $2$ which we have shown to be false in the previous subsection (via an argument similar to the above one).
\end{proof}
\end{lem}
Hence, we can conclude from the lemma that in order to make the action (with generic orbit of dimension $2$) equidimensional it suffices to blow up compact surfaces consisting of fixed points. Let us now show that this process preserves de Rham cohomology.
\begin{lem} Let $M$ be a compact connected smooth $5$-manifold with a (smooth) $\mathbb{S}^3$-action with generic orbit of dimension $2$ which contains a surface consisting of fixed points. Let us denote by $\hat{M}$ the blow up of that surface in $M$. Then $H^{\bullet}_{dR}(M)\cong H^{\bullet}_{dR}(\hat{M})$.
\begin{proof} The proof is verbatim the same as in the analogous lemma in the previous subsection (with $\mathbb{S}^1$ replaced by the surface of fixed points).
\end{proof}
\end{lem}
To finish the proof it suffices to consider $\hat{M}$ (the blow up of all surfaces of fixed points in $M$) by the above lemma. Note that its second cohomology contains $E^{2,0}_2(\hat{M})$ as a subspace (since all the coboundary operators entering and exiting this position are zero on all subsequent pages). The hypersurface on which the product vanishes corresponds to $E^{2,0}_2(\hat{M})$. To see this note that $E^{1,1}_2(\hat{M})=0$ and $E^{0,2}_2(\hat{M})$ is either $0$ or $\mathbb{R}$. If it is zero or this term vanishes before the sequence collapses then $E^{2,0}_2(\hat{M})$ has codimension zero and the statement would be true by taking a hypersurface within it. Otherwise the subpsace corresponding to $E^{2,0}_2(\hat{M})$ in the second cohomology is precisely the required hypersurface. To see that the product indeed vanishes when restricted to this subspace it suffices to note that multiplying two elements from it will result in an element from $E^{4,0}_2(\hat{M})\cong 0$. Consequently, the product on the second page has to vanish and since this is a first quadrant spectral sequence and the result of this product is on the $0$-th row the same is true for the product of the corresponding elements in the second cohomology.
\newline\indent To treat the case when the generic orbit has dimension $3$ we will need to systematically study the various possibilities for the strata of orbits of smaller dimension and show that the process of blowing up only changes the second real homology in one particular case in which it reduces it by a class represented by a sphere. We will partially work with $H_{dR}^{\bullet}(M)$ with the understanding that by the Universal Coefficients Theorem and the de Rham isomorphism it is isomorphic to $H_{\bullet}(M,\mathbb{R})$. Let us start by analyzing how blowing up the isolated fixed points of the action influences the cohomology of the manifold
\begin{lem} Let $M$ be a compact connected smooth $5$-manifold with a (smooth) $\mathbb{S}^3$-action with generic orbit of dimension $3$ which contains an isolated fixed point. Let us denote by $\hat{M}$ the blow up of that fixed point in $M$. Then $H^{\bullet}_{dR}(M)\cong H^{\bullet}_{dR}(\hat{M})$.
\begin{proof}  Let us again consider the blow down map $f:\hat{M}\to M$ identifying the copy of $\mathbb{R}P^4$ to a point. For this map let us write down the Mayer-Vietoris ladder diagram:
$$\begin{tikzcd}
... \arrow[r]& H^{k-1}_{dR}(U)\oplus H^{k-1}_{dR}(V)\arrow[r]\arrow[d]& H^{k-1}_{dR}(U\cap V)\arrow[r]\arrow[d]& H^{k}_{dR}(M)\arrow[r]\arrow[d]&H^{k}_{dR}(U)\oplus H^{k}_{dR}(V)\arrow[r]\arrow[d]& ...\\
... \arrow[r]& H^{k-1}_{dR}(\hat{U})\oplus H^{k-1}_{dR}(V)\arrow[r]& H^{k-1}_{dR}(\hat{U}\cap V)\arrow[r]& H^{k}_{dR}(\hat{M})\arrow[r]&H^{k}_{dR}(\hat{U})\oplus H^{k}_{dR}(V)\arrow[r]& ...
\end{tikzcd}$$
where $V$ is just $\hat{M}\backslash\mathbb{R}P^4$ which can be identified with a subset of $M$ as well, $U$ is a small ball around the isolated fixed point and $\hat{U}=f^{-1}(U)$. We are going to use the five lemma to show that $f$ induces an isomorphism in cohomology. Let us note that the map on $V$ is the identity and hence induces an isomorphism. Moreover, $f$ restricts to a diffeomorphism from $\hat{U}\cap V$ to $U\cap V$. Hence, it suffices to show that $f|_{\hat{U}}:\hat{U}\to U$ induces an isomorphism in cohomology. By homotopy invariance it suffices to prove that this is true for $f|_{\mathbb{R}P^4}$. This however is obvious since this map simply collapses $\mathbb{R}P^4$ to a point and $\mathbb{R}P^4$ has the same real cohomology as a point and consequently it induces an isomorphism in cohomology finishing the proof by the five lemma.
\end{proof}
\end{lem}
Let us further note that if in $H_2(\hat{M},\mathbb{R})$ there is a hypersurface generated by embedded spheres then the same is true for $H_2(M,\mathbb{R})$. To see this note that the image of any class represented by an embedded sphere in $\hat{M}$ is again an image of the fundamental class of a sphere through the composition. Furthermore, since $M$ is of dimension $5$ the image of this composition can be changed by a homotopy to an embedding. Hence, it suffices to prove the Theorem for $M$ with all fixed points blown up.
\newline\indent By the above consideration we assume without loss of generality that $M$ has no isolated fixed points. This implies that the fixed point strata consists of circles (by a dimension argument as before). We need to show that there are no $2$-dimensional orbits arbitrarily close to such a circle.
\begin{lem} Let $M$ be a compact connected smooth $5$-manifold with a (smooth) $\mathbb{S}^3$-action with generic orbit of dimension $3$ and without isolated fixed points. Then there exists an open neighbourhood of the fixed point strata which is disjoint from the set of points belonging to orbits of dimension $2$.
\begin{proof} Let us chose a fixed point $x$ on such a circle and an invariant metric on $M$. Then the action preserves the geodesic $4$-disc perpendicular to $\mathbb{S}^1$ through $x$ and its boundary $\mathbb{S}^3$. Moreover, if the radius of the disc is small enough then there are no fixed points of the action restricted to this $\mathbb{S}^3$. Hence, since we have already proved that there is no action of $\mathbb{S}^3$ on $\mathbb{S}^3$ with all orbits of dimension $2$ the generic orbit has to be of dimension $3$. Hence, there are no $2$-dimensional orbits close to such a circle.
\end{proof}
\end{lem}
This allows us to consider the fixed point strata and the two strata corresponding to $2$-dimensional orbits seperately. Let us now show that blowing up the circles of fixed points does not change $H_2(M,\mathbb{R})$.
\begin{lem} Let $M$ be a compact connected smooth $5$-manifold with a (smooth) $\mathbb{S}^3$-action with generic orbit of dimension $3$ which contains a circle of fixed points. Let us denote by $\hat{M}$ the blow up of that circle in $M$. Then $H^{k}_{dR}(M)\cong H^{k}_{dR}(\hat{M})$ for $k\leq 2$.
\begin{proof}  Let us consider the smooth map $f:\hat{M}\to M$ (blow down) defined using the second description of the blow up process as follows. On $V:=M\backslash  \overline{T}$ (where $T$ is the tubular saturated neighbourhood of the circle) it is simply the identity. On $\hat{U}:=P_{NS}\times_{O(4)}\hat{\mathbb{R}}^4\to P_{NS}\times_{O(4)}\mathbb{R}^4=:U$ it is given by identity on the first coordinate and collapsing $\mathbb{R}P^3$ to the origin in $\mathbb{R}^4$ on the second coordinate. For this map let us write down the Mayer-Vietoris ladder diagram:
$$\begin{tikzcd}
... \arrow[r]& H^{k-1}_{dR}(U)\oplus H^{k-1}_{dR}(V)\arrow[r]\arrow[d]& H^{k-1}_{dR}(U\cap V)\arrow[r]\arrow[d]& H^{k}_{dR}(M)\arrow[r]\arrow[d]&H^{k}_{dR}(U)\oplus H^{k}_{dR}(V)\arrow[r]\arrow[d]& ...\\
... \arrow[r]& H^{k-1}_{dR}(\hat{U})\oplus H^{k-1}_{dR}(V)\arrow[r]& H^{k-1}_{dR}(\hat{U}\cap V)\arrow[r]& H^{k}_{dR}(\hat{M})\arrow[r]&H^{k}_{dR}(\hat{U})\oplus H^{k}_{dR}(V)\arrow[r]& ...
\end{tikzcd}$$
We are going to use the five lemma to show that $f$ induces an isomorphism in cohomology up to degree $2$. For this let us note that the map on $V$ is the identity and hence induces an isomorphism. Moreover, $f$ restricts to a diffeomorphism from $\hat{U}\cap V$ to $U\cap V$. Hence, it suffices to show that $f|_{\hat{U}}:\hat{U}\to U$ induces an isomorphism in cohomology for degree up to $2$. By homotopy invariance it suffices to prove that this is true for $f|_{P_{NS}\times_{O(4)}\hat{\mathbb{R}P}^3}:P_{NS}\times_{O(4)}\mathbb{R}P^3\to\mathbb{S}^1$. This however follows since from the Serre spectral sequence $P_{NS}\times_{O(4)}\mathbb{R}P^3$ has the same cohomology up to degree $2$ as $\mathbb{S}^1$ while this map is simply the projection in this bundle. To see this note that this sequence degenerates at the second page (and as before the action of the fundamental group of the base on the cohomology of the fiber is trivial) with the bottom row (which contains the only non-zero elements contributing to cohomology up to degree $2$) corresponding to the cohomology of the base and hence the resulting cohomology classes induced by that row are simply the pullback of the classes from the circle. Consequently, this map induces an isomorphism in cohomology up to degree $2$ finishing the proof by the five lemma.
\end{proof}
\end{lem}
Same as before we can conclude that if this blow up has the desired hypersurface in the second homology then so does the initial manifold $M$ and we can assume without loss of generality that the action on $M$ has no fixed points. Let us now exclude the presence of isolated two dimensional orbits.
\begin{lem} Let $M$ be a compact connected smooth $5$-manifold with a (smooth) $\mathbb{S}^3$-action with generic orbit of dimension $3$ without fixed points. Then any component of the $\mathbb{S}^2$ strata is either a $\mathbb{S}^2$-bundle over a circle or equivariantly diffeomorphic to $\mathbb{S}^2\times (0,1)$ with the orbits at the ends having type $\mathbb{R}P^2$. The components of the $\mathbb{R}P^2$ strata are either isolated $\mathbb{R}P^2$ type orbits at the ends of $\mathbb{S}^2\times (0,1)$ or $\mathbb{R}P^2$ bundles over $\mathbb{S}^1$.
\begin{proof} For an orbit of type $\mathbb{S}^2$ the corresponding isotropy group is $\mathbb{S}^1$ and the normal bundle is $3$-dimensional. Describing the action of $\mathbb{S}^1$ on the fiber through the chosen point amounts to a non-trivial Lie group homomorphism from $\mathbb{S}^1$ to $O(3)$. Hence, up to a covering of $\mathbb{S}^1$ this action is just the standard action of a $1$-parameter subgroup of $O(3)$. But these are just rotations around the axis and hence there is a $1$-dimensional subspace fixed by $\mathbb{S}^1$. Hence, by the Slice Theorem this strata is locally just a $\mathbb{S}^2\times (0,1)$. These have to either glue to a $\mathbb{S}^2$-bundle over a circle or there has to be a further quotient of $\mathbb{S}^2$ at the end of $\mathbb{S}^2\times (0,1)$. In the later case the only possibility is $\mathbb{R}P^2$ (since there are no fixed points).
\newline\indent Let us now assume that there is an orbit of type $\mathbb{R}P^2$. Then the isotropy group $H$ consists of two circles and we study possible non-trivial homomorphisms of this group into $O(3)$. The image of the identity component is again just a $1$-parameter subgroup whereas the image of the second circle has to be one of the disjoint circles in the normalizer of this $1$-parameter subgroup (possibly the same one). These circles correspond to rotations along the specified axis composed with the antipodal map or the symmetry with respect to some plane containing the axis of rotation (possibly neither or both). Regardless of which of these is the image of the other circle it also preserves the $1$-dimensional space corresponding to the axis of rotation. Hence, if the image of the other circle isn't one of the circles corresponding to composition with the antipodal map then the other orbits corresponding to points on this line have to be of the same type as the central orbit. Otherwise, they will be copies of $\mathbb{S}^2$.
\end{proof}
\end{lem}
\begin{rem} In fact, the case considered in the last sentence of the proof also shows that near such an orbit of type $\mathbb{R}P^2$ at the end of an $\mathbb{S}^2\times (0,1)$ there can be no other $2$-dimensional orbits. Hence, in this case the closure of $\mathbb{S}^2\times (0,1)$ is already a manifold and hence we can immediately blow it up.
\end{rem}
Let us consider the remaining options for components of the set of two dimensional orbits. Let us start by treating the case of blowing up the closure of $\mathbb{S}^2\times (0,1)$:
\begin{lem} Let $M$ be a compact connected smooth $5$-manifold with a (smooth) $\mathbb{S}^3$-action with generic orbit of dimension $3$ (and no fixed points) with component of the $\mathbb{S}^2$ strata equivariantly diffeomorphic to $\mathbb{S}^2\times (0,1)$. Let us denote by $\hat{M}$ the blow up of the closure of such a component in $M$. Then $H^{k}_{dR}(M)\cong H^{k}_{dR}(\hat{M})$ for $k\leq 2$.
\begin{proof}  Let us consider the smooth map $f:\hat{M}\to M$ (blow down) defined using the second description of the blow up process as follows. On $V:=M\backslash  \overline{T}$ (where $T$ is the tubular saturated neighbourhood of the closure of this component of the strata) it is simply the identity. On $\hat{U}:=P_{NS}\times_{O(2)}\hat{\mathbb{R}}^2\to P_{NS}\times_{O(2)}\mathbb{R}^2=:U$ it is given by identity on the first coordinate and collapsing $\mathbb{R}P^1$ to the origin in $\mathbb{R}^2$ on the second coordinate. For this map let us write down the Mayer-Vietoris ladder diagram:
$$\begin{tikzcd}
... \arrow[r]& H^{k-1}_{dR}(U)\oplus H^{k-1}_{dR}(V)\arrow[r]\arrow[d]& H^{k-1}_{dR}(U\cap V)\arrow[r]\arrow[d]& H^{k}_{dR}(M)\arrow[r]\arrow[d]&H^{k}_{dR}(U)\oplus H^{k}_{dR}(V)\arrow[r]\arrow[d]& ...\\
... \arrow[r]& H^{k-1}_{dR}(\hat{U})\oplus H^{k-1}_{dR}(V)\arrow[r]& H^{k-1}_{dR}(\hat{U}\cap V)\arrow[r]& H^{k}_{dR}(\hat{M})\arrow[r]&H^{k}_{dR}(\hat{U})\oplus H^{k}_{dR}(V)\arrow[r]& ...
\end{tikzcd}$$
We are going to use the five lemma to show that $f$ induces an isomorphism in cohomology up to degree $2$. For this let us note that the map on $V$ is the identity and hence induces an isomorphism. Moreover, $f$ restricts to a diffeomorphism from $\hat{U}\cap V$ to $U\cap V$. Hence, it suffices to show that $f|_{\hat{U}}:\hat{U}\to U$ induces an isomorphism in cohomology for degree up to $2$. By homotopy invariance it suffices to prove that this is true for $f|_{P_{NS}\times_{O(2)}\mathbb{R}P^1}:P_{NS}\times_{O(2)}\mathbb{R}P^1\to S$ (where $S$ denotes the closure of $\mathbb{S}^2\times (0,1)$). Computing the de Rham cohomology of $S$ via a Mayer-Vietoris sequence (with respect to the neighbourhoods of the two copies of $\mathbb{R}P^2$) we readilly get that it is isomorphic to that of $\mathbb{S}^3$. On the other hand, $P_{NS}\times_{O(2)}\mathbb{R}P^1$ is a closed $4$-manifold with all the orbits of dimension $3$. Moreover, to every orbit in $S$ corresponds a single orbit in $P_{NS}\times_{O(2)}\mathbb{R}P^1$ for dimensional reasons. Hence, $P_{NS}\times_{O(2)}\mathbb{R}P^1$ consists of $\mathbb{S}^3\slash\Gamma\times (0,1)$ with some $3$-dimensional orbits at the ends. By cross checking this with the analysis of cases in the previous subsection (for $4$-manifolds with generic orbit of dimension $3$) we see that manifolds fitting this description all have the same de Rham cohomologyy as $\mathbb{S}^3$ up to degree $2$. Hence, the map restricted to $\hat{U}$ induces an isomorphism in cohomology up to degree $2$ and consequently by the five lemma the same is true for the entire map induced by $f$ finishing the proof.
\end{proof}
\end{lem}
Performing the blow up in the case of bundles with fiber $\mathbb{R}P^2$ will again not change the second homology:
\begin{lem} Let $M$ be a compact connected smooth $5$-manifold with a (smooth) $\mathbb{S}^3$-action with generic orbit of dimension $3$ (and no fixed points) with component of the $\mathbb{R}P^2$ strata equal to an $\mathbb{R}P^2$ bundle over a circle. Let us denote by $\hat{M}$ the blow up of that bundle in $M$. Then $H^{k}_{dR}(M)\cong H^{k}_{dR}(\hat{M})$ for $k\leq 2$.
\begin{proof}  Let us consider the smooth map $f:\hat{M}\to M$ (blow down) defined using the second description of the blow up process as follows. On $V:=M\backslash  \overline{T}$ (where $T$ is the tubular saturated neighbourhood of this bundle) it is simply the identity. On $\hat{U}:=P_{NS}\times_{O(2)}\hat{\mathbb{R}}^2\to P_{NS}\times_{O(2)}\mathbb{R}^2=:U$ it is given by identity on the first coordinate and collapsing $\mathbb{R}P^1$ to the origin in $\mathbb{R}^2$ on the second coordinate. For this map let us write down the Mayer-Vietoris ladder diagram:
$$\begin{tikzcd}
... \arrow[r]& H^{k-1}_{dR}(U)\oplus H^{k-1}_{dR}(V)\arrow[r]\arrow[d]& H^{k-1}_{dR}(U\cap V)\arrow[r]\arrow[d]& H^{k}_{dR}(M)\arrow[r]\arrow[d]&H^{k}_{dR}(U)\oplus H^{k}_{dR}(V)\arrow[r]\arrow[d]& ...\\
... \arrow[r]& H^{k-1}_{dR}(\hat{U})\oplus H^{k-1}_{dR}(V)\arrow[r]& H^{k-1}_{dR}(\hat{U}\cap V)\arrow[r]& H^{k}_{dR}(\hat{M})\arrow[r]&H^{k}_{dR}(\hat{U})\oplus H^{k}_{dR}(V)\arrow[r]& ...
\end{tikzcd}$$
We are going to use the five lemma to show that $f$ induces an isomorphism in cohomology up to degree $2$. For this let us note that the map on $V$ is the identity and hence induces an isomorphism. Moreover, $f$ restricts to a diffeomorphism from $\hat{U}\cap V$ to $U\cap V$. Hence, it suffices to show that $f|_{\hat{U}}:\hat{U}\to U$ induces an isomorphism in cohomology for degree up to $2$. By homotopy invariance it suffices to prove that this is true for $f|_{P_{NS}\times_{O(2)}\mathbb{R}P^1}:P_{NS}\times_{O(2)}\mathbb{R}P^1\to \mathbb{R}P^2\times\mathbb{S}^1$ (since the bundle we blow up has to be glued by equivariant diffeomorphisms of $\mathbb{R}P^2$ and hence is just the product). Note that when we restrict the above $\mathbb{R}P^1\cong\mathbb{S}^1$ bundle to any $\mathbb{R}P^2\times\{x\}$ we have to get a single $3$ dimensional orbit and since the bundles over different copies are isomorphic the resulting orbits are equivariantly diffeomorphic. Hence, as in the previous subsection we have an $\mathbb{S}^3\slash\Gamma$ bundle over $\mathbb{S}^1$ with trivial action of the fundamental group on the fibers. Hence, its cohomology up to degree two agrees with that of $\mathbb{R}P^2\times\mathbb{S}^1$. Moreover, since in both cases the zeroth and first cohomology is generated by pullbacks of forms from the base circle, which is preserved by $f|_{P_{NS}\times_{O(2)}\mathbb{R}P^1}$, we can see that this map indeed induces an isomorphism in cohomology up to degree $2$. This finishes the proof by the five lemma.
\end{proof}
\end{lem}
Finally, in the case of $\mathbb{S}^2$ bundles over $\mathbb{S}^1$ we will show that the only change in second homology that can occur is lowering the dimension by $1$ which corresponds to the space generated by the fiber of this bundle in $M$. For this let $\hat{M}$ denote the blow up of this components of the strata and let us one final time consider the smooth map $f:\hat{M}\to M$ (blow down) defined using the second description of the blow up process as follows. On $V:=M\backslash  \overline{T}$ (where $T$ is the tubular saturated neighbourhood of this bundle) it is simply the identity. On $\hat{U}:=P_{NS}\times_{O(2)}\hat{\mathbb{R}}^2\to P_{NS}\times_{O(2)}\mathbb{R}^2=:U$ it is given by identity on the first coordinate and collapsing $\mathbb{R}P^1$ to the origin in $\mathbb{R}^2$ on the second coordinate. Again let us write its Mayer-Vietoris ladder diagram:
$$\begin{tikzcd}
... & H_{k-1}(\hat{U},\mathbb{R})\oplus H_{k-1}(V,\mathbb{R})\arrow[l]\arrow[d]& H_{k-1}(\hat{U}\cap V,\mathbb{R})\arrow[l]\arrow[d]& H_{k}(\hat{M},\mathbb{R})\arrow[l]\arrow[d]& ...\arrow[l]\\
... & H_{k-1}(U,\mathbb{R})\oplus H_{k-1}(V,\mathbb{R})\arrow[l]& H_{k-1}(U\cap V,\mathbb{R})\arrow[l]& H_{k}(M,\mathbb{R})\arrow[l]& ...\arrow[l]
\end{tikzcd}$$
Since there are two possible $\mathbb{S}^2$ bundles over $\mathbb{S}^1$ admitting a $\mathbb{S}^3$ action (the trivial one and one glued by the antipodal map which is the only non-trivial $\mathbb{S}^3$-equivariant diffeomorphism of $\mathbb{S}^2$). One can compute from the Serre spectral sequence (or our spectral sequence) that the de Rham cohomology of the non-trivial case is again the same as that of a circle and hence by a verbatim the same argument as in the previous lemma second cohomology remains unchanged. For the trivial case we can still conclude that the map on $V$ is the identity and hence induces an isomorphism and $f$ restricts to a diffeomorphism from $\hat{U}\cap V$ to $U\cap V$. Hence, it suffices to show that $f|_{\hat{U}}:\hat{U}\to U$ induces an isomorphism in cohomology for degree up to $1$ and then study its behaviour in degree $2$. By homotopy invariance it suffices to prove that this is true for $f|_{P_{NS}\times_{O(2)}\mathbb{R}P^1}:P_{NS}\times_{O(2)}\mathbb{R}P^1\to \mathbb{S}^2\times\mathbb{S}^1$. Note that when we restrict the above $\mathbb{R}P^1\cong\mathbb{S}^1$ bundle to any $\mathbb{S}^2\times\{x\}$ we have to get a single $3$ dimensional orbit and since the bundles over different copies are isomorphic the resulting orbits are equivariantly diffeomorphic. Hence, as in the proof of the previous lemma we have an $\mathbb{S}^3\slash\Gamma$ bundle over $\mathbb{S}^1$ with trivial action of the fundamental group on the fibers. Hence, its cohomology up to degree one agrees with that of $\mathbb{S}^2\times\mathbb{S}^1$. Moreover, since in both cases the zeroth and first cohomology is generated by pullbacks of forms from the base circle, which is preserved by $f|_{P_{NS}\times_{O(2)}\mathbb{R}P^1}$, we can see that this map indeed induces an isomorphism in cohomology up to degree $1$ and hence the same is true with real homology up to degree $1$ (by the Universal Coefficients Theorem). In degree, $2$ however this map is the $0$ map. Hence, after blowing up the second homology of $U$ change from $\mathbb{R}$ to $0$. Since up to degree one these long exact sequences are isomorphic this may result in one of two things in the cohomology of $\hat{M}$:
\begin{enumerate}
\item If the image of the class of the fiber $\mathbb{S}^2$ (generating $H_2(U,\mathbb{R})$) in $M$ (under the inclusion of $U$) is zero then the blow up results in increasing the dimension of the third homology by one.
\item If the image of this class is non-zero then this class is represented by an embedded sphere in $M$ and the blow up will result in decreasing the dimension of the second homology by one with this class generating the complement of the image through $f$ of $H_2(\hat{M},\mathbb{R})$ in $H_{2}(M,\mathbb{R})$.
\end{enumerate}
Now the proof of the second case can be finished by using our spectral sequence. Namely, for an equidimensional action of $\mathbb{S}^3$ with generic orbit of dimension $3$, one can note that in our spectral sequence both $E^{0,2}_2$ and $E^{1,1}_2$ vanish whereas $E^{2,0}_2\cong E^{2,0}_{\infty}$ is either $0$ or $\mathbb{R}$ (since the foliation by orbits is Riemannian). Hence, the second cohomology (and consequently the second real homology by the Universal Coefficients Theorem) is either $0$ or $\mathbb{R}$. This leads us to the conclusion that the process of blowing up an arbitrary action on a given manifold has to reduce the second homology so that its dimension is at most $1$. But as we have shown above all of the reduced classes had to be represented by embedded spheres and so they form the desired hypersurface.
\subsection{Comments on the relevance of the results of this section.} In this subsection we wish to provide some context for the result from the previous subsection. This last section is heavilly inspired by \cite{Spin} on which we base the following exposition and highly recommend for further information.
\begin{prop}[Proposition 3.1, Chapter IV of \cite{Spin}]\label{NAA} Let $G$ be a compact connected Lie group of dimension greater then zero. Then $G$ contains $\mathbb{S}^1$ as a subgroup. In fact, such subgroups are dense in $G$. Moreover, if $G$ is non-abelian, then it contains a closed Lie subgroup isomorphic to a quotient of $\mathbb{S}^3$ by some finite normal group.
\end{prop}
This informs the following result:
\begin{prop}\label{NAA2} A manifold $M$ admits an effective action of some compact connected non-abelian Lie group if and only if it admits a non-trivial (on every connected component) action of $\mathbb{S}^3$.
\begin{proof} Let us note that non-trivial normal subgroups of $\mathbb{S}^3$ are finite and hence an effective action of a quotient of $\mathbb{S}^3$ by some finite normal group is equivalent to a non-trivial action of $\mathbb{S}^3$. Hence, since such a quotient of $\mathbb{S}^3$ is a compact connected non-abelian Lie group one implication is obvious. On the other hand, by the previous Theorem if there is some compact connected non-abelian Lie group acting effectively on $M$ then it has a subgroup isomorphic to a quotient of $\mathbb{S}^3$ by some finite normal group. Hence, restricting the action to this subgroup gives us the desired effective action of a quotient of $\mathbb{S}^3$ by some finite normal subgroup (or equivalently a non-trivial action of $\mathbb{S}^3$).
\end{proof}
\end{prop}
Hence, as we have already seen in the previous subsection, the study of $\mathbb{S}^3$-actions can provide very general information about the existence of actions of arbitrary groups on the given manifold. Some further (non-obvious) justification for studying such actions or lack thereof is provided by the following important result linking it to the classical and thoroughly studied problem (see e.g. \cite{Spin}) of the existence of a metric with positive scalar curvature on a given smooth manifold $M$:
\begin{tw}[Theorem 8.14, Chapter II of \cite{Spin}] Let $M$ be a compact manifold which admits an effective (smooth) action by a compact, connected, non-abelian Lie group. Then $M$ admits a metric of positive scalar curvature.
\end{tw}
As we have shown earlier and is pointed out in the original paper (see \cite{LY}) the assumption on the existence of a Lie group acting effectively on $M$ can be replaced by the existence of a non-trivial action of $\mathbb{S}^3$ giving us the following equivalent version of this Theorem:
\begin{cor} Let $M$ be a compact manifold which admits a non-trivial (on every connected component) action of $\mathbb{S}^3$. Then $M$ admits a metric of positive scalar curvature.
\end{cor}
It is remarked in \cite{Spin} that these results provide an obstruction to the existence of Lie group actions on manifolds (either in term of the existence of metrics with positive scalar curvature directly or more practically in terms of the $\hat{A}$-genus via the celebrated Atiyah-Singer Index Theorem). However, in principal one can try to construct a non-trivial action of $\mathbb{S}^3$ on a given manifold in order to show (via the above results) that it admits a metric of positive scalar curvature. For such an approach it would be convenient to have a tool allowing us to exclude manifolds which cannot admit any such group actions. While it is obvious from the nature of the sequence and the above results that our sequence can serve this purpose on a case by case basis, the results presented in this section demonstrates that it can also exclude large classes of manifolds all at once.
\begin{rem} The value of the approach proposed above stems from it not relying on the existence of spin structures (unlike criterions based on the Atiyah-Singer Index Theorem). Hence, we believe it is worth further development.
\end{rem}

\bigskip
\Small{\textbf{Statements and Declarations} The authors declare no competing interests.}
\bigskip

\Small{\textbf{Availability of data and material} Data sharing not applicable to this article as no datasets were generated or analysed during the current study.}


\begin{thebibliography}{99}
\bibitem{ss1} J.A. \'Alvarez L\'opez;
\emph{A finiteness theorem for the spectral sequence of a Riemannian foliation}. Illinois J. Math. 33(1): 79-92 (1989).
\bibitem{ss2}  J.A. \'Alvarez L\'opez;
\emph{A Decomposition Theorem for the Spectral Sequence of Lie Foliations}. T. Am. Math. Soc. Vol. 329, No. 1, 173-184 (1992).
\bibitem{ss3} J.A. \'Alvarez L\'opez; Y.A. Kordyukov;
\emph{Adiabatic limits and spectral sequence for Riemannian foliations}. GAFA, Geom. funct. anal. 10, 977-1027 (2000).
\bibitem{BR-A} M. Benameur; A. Rey-Alcantara;
\emph{La signature basique est un invariant d’homotopie feuillet´ee}. C. R. Math. Acad. Sci. Paris 349, no. 13-14, 787–791 (2011).
\bibitem{Spin} H. B. Lawson, Jr.; M.-L. Michelsohn;
\emph{Spin Geometry}. Princeton Mathematical Series, vol. 38, Princeton University Press (1989). 
\bibitem{LY} H. B. Lawson, Jr.; S.T. Yau;
\emph{Scalar curvature, non-abelian group actions, and the degree of symmetry of exotic spheres}. Comm. Math. Helv. 49, 232-244 (1974).
\bibitem{ttt} R. Bott; L.W. Tu;
\emph{Differemtial forms in algebraic topology}. Springer New York (1982).
\bibitem{BG} C.P. Boyer; K. Galicki;
\emph{Sasakian Geometry}. Oxford Mathematical Monographs. Oxford University Press (2007).
\bibitem{B} G.E. Bredon;
\emph{Introduction to Compact Transformation Groups}. Pure Appl. Math. vol. 42, Academic Press N.Y (1972).
\bibitem{E1} A. El Kacimi-Alaoui;
\emph{Op\'{e}rateurs transversalement elliptiques sur un feuilletage riemannien et applications}. Compositio Mathematica, 73, 57-106 (1990).
\bibitem{HR} G. Habib; K. Richardson;
\emph{Homotopy invariance of cohomology and signature of a Riemannian foliation}. Math. Z. 293, 579–595 (2019).
\bibitem{Hoch1} G. Hochschild;
\emph{The Automorphism Group of a Lie Group} T. Am. Math. Soc.. 72, 2, pp. 209-216 (1952).
\bibitem{Hoch2} G. Hochschild; J.-P. Serre;
\emph{Cohomology of Lie algebras} Ann. of Math. 57, 3, pp. 591-603 (1953).
\bibitem{Str} J.N. Mather;
\emph{Stratification and Mappings}. Dynamical Systems, Academic Press N.Y. 9, 12, 195-232 (1973).
\bibitem{MP} P. Melvin; J. Parker;
\emph{4-manifolds with large symmetry groups}. Topology 25, no. 1, 71–83 (1986).
\bibitem{M1}J.W. Milnor;
\emph{Curvatures of left invariant metrics on lie groups}. Advances in Mathematics, 21:293–329 (1976).
\bibitem{My2} P. Ra\'{z}ny;
\emph{Invariance of basic numbers under deformations of Sasakian manifolds}. ANN MAT PUR APPL ,200, pages 1451–1468 (2021).
\bibitem{My3} P. Raźny;
\emph{Cohomology of manifolds with structure group $U(n)\times O(s)$}. Geom Dedicata 217, 58 (2023).
\bibitem{My4} P. Raźny;
\emph{A Spectral Sequence for Locally Free Isometric Lie Group Actions }.  Transformation Groups (2024)  https://doi.org/10.1007/s00031-024-09855-2 .
\bibitem{RPSA} J.I. Royo Prieto; M. Saralegui-Aranguren;
\emph{The Gysin sequence for $\mathbb{S}^3$-actions on manifolds}. PUBL MATH-DEBRECEN  83, 3, 275-289 (2013). 
\bibitem{W} J.I. Royo Prieto; M. Saralegui-Aranguren; R. Wolak;
\emph{Hard Lefschetz Property for Isometric Flows}. Transformation Groups (2022) https://doi.org/10.1007/s00031-022-09744-6.
\bibitem{ss4} K.S. Sarkaria;
\emph{A finiteness theorem for foliated manifolds}. J. Math. Soc. Japan 30, 687-696 (1978).
\bibitem{T} A. Tralle; J. Oprea;
\emph{Symplectic Manifolds with no K\"{a}hler Structure}. Lect. Notes
Math. 1661, Springer, Berlin (1997).
\end{thebibliography}
\end{document}